\documentclass[reqno,11pt]{article}
\usepackage{mathrsfs}
\usepackage{mathtools}
\usepackage{amssymb}
\usepackage{color}
\usepackage{amsthm}
\usepackage{amsmath}
\usepackage{amsfonts}
\usepackage{enumerate}
\usepackage{bm}

\usepackage{graphicx}

\numberwithin{equation}{section}

\newtheorem{theorem}{Theorem}[section]
\newtheorem{lemma}[theorem]{Lemma}
\newtheorem{corollary}[theorem]{Corollary}
\newtheorem{proposition}[theorem]{Proposition}

\theoremstyle{definition}
\newtheorem{definition}[theorem]{Definition}
\newtheorem{assumption}[theorem]{Assumption}
\newtheorem{example}[theorem]{Example}

\theoremstyle{remark}
\newtheorem{remark}[theorem]{Remark}

\begin{document}
\title{Decay properties of Axially Symmetric D-solutions to the Steady Incompressible Magnetohydrodynamic Equations}
\author{Shangkun Weng\thanks{School of mathematics and statistics, Wuhan University, Wuhan, Hubei Province, 430072, People's Republic of China. Email: skweng@whu.edu.cn}\and Yan Zhou\thanks{School of mathematics and statistics, Wuhan University, Wuhan, Hubei Province, 430072, People's Republic of China. Email: yanz0013@whu.edu.cn}}
\date{}
\maketitle

\def\be{\begin{eqnarray}}
\def\ee{\end{eqnarray}}
\def\ba{\begin{aligned}}
\def\ea{\end{aligned}}
\def\bay{\begin{array}}
\def\eay{\end{array}}
\def\bca{\begin{cases}}
\def\eca{\end{cases}}
\def\p{\partial}
\def\hphi{\hat{\phi}}
\def\bphi{\bar{\phi}}
\def\no{\nonumber}
\def\eps{\epsilon}
\def\de{\delta}
\def\De{\Delta}
\def\om{\omega}
\def\Om{\Omega}
\def\f{\frac}
\def\th{\theta}
\def\vth{\vartheta}
\def\la{\lambda}
\def\lab{\label}
\def\b{\bigg}
\def\var{\varphi}
\def\na{\nabla}
\def\ka{\kappa}
\def\al{\alpha}
\def\La{\Lambda}
\def\ga{\gamma}
\def\Ga{\Gamma}
\def\ti{\tilde}
\def\wti{\widetilde}
\def\wh{\widehat}
\def\ol{\overline}
\def\ul{\underline}
\def\Th{\Theta}
\def\si{\sigma}
\def\Si{\Sigma}
\def\oo{\infty}
\def\q{\quad}
\def\z{\zeta}
\def\co{\coloneqq}
\def\eqq{\eqqcolon}
\def\di{\displaystyle}
\def\bt{\begin{theorem}}
\def\et{\end{theorem}}
\def\bc{\begin{corollary}}
\def\ec{\end{corollary}}
\def\bl{\begin{lemma}}
\def\el{\end{lemma}}
\def\bp{\begin{proposition}}
\def\ep{\end{proposition}}
\def\br{\begin{remark}}
\def\er{\end{remark}}
\def\bd{\begin{definition}}
\def\ed{\end{definition}}
\def\bpf{\begin{proof}}
\def\epf{\end{proof}}
\def\bex{\begin{example}}
\def\eex{\end{example}}
\def\bq{\begin{question}}
\def\eq{\end{question}}
\def\bas{\begin{assumption}}
\def\eas{\end{assumption}}
\def\ber{\begin{exercise}}
\def\eer{\end{exercise}}
\def\mb{\mathbb}
\def\mbR{\mb{R}}
\def\mbZ{\mb{Z}}
\def\mc{\mathcal}
\def\mcS{\mc{S}}
\def\ms{\mathscr}
\def\lan{\langle}
\def\ran{\rangle}
\def\lb{\llbracket}
\def\rb{\rrbracket}
\def\fr#1#2{{\frac{#1}{#2}}}
\def\dfr#1#2{{\dfrac{#1}{#2}}}
\def\u{{\textbf u}}
\def\v{{\textbf v}}
\def\h{\textbf{h}}
\def\w{{\textbf w}}
\def\d{{\textbf d}}
\def \ff{\textbf{f}}
\def\g{\textbf{g}}
\def\nn{{\textbf n}}
\def\x{{\textbf x}}
\def\e{{\textbf e}}
\def\D{{\textbf D}}
\def\U{{\textbf U}}
\def\M{{\textbf M}}
\def\F{{\mathcal F}}
\def\I{{\mathcal I}}
\def\W{{\mathcal W}}
\def\div{{\rm div\,}}
\def\curl{{\rm curl\,}}
\def\R{{\mathbb R}}
\def\FF{{\textbf F}}
\def\A{{\textbf A}}
\def\R{{\textbf R}}
\def\r{{\textbf r}}
\def\={\coloneqq}

\begin{abstract}
  In this paper, we investigate the decay properties of axially symmetric solutions to the steady incompressible magnetohydrodynamic equations in $\mbR^3$ with finite Dirichlet integral. We first derive the decay rates of general D-solutions to the axisymmetric MHD equations. In the special case where the magnetic field only has the swirl component ${\bf h}(r,z)= h_{\th}(r,z) {\bf e_{\th}}$, we obtain better decay rates. The last result examines the decay rates along the axis $Oz$ also within the special class of D-solutions with only swirl magnetic field. The main tool in this paper is the combination of the scaling argument, the \textit{Brezis-Gallouet inequality} and the weighted energy estimate.
\end{abstract}

\begin{center}
\begin{minipage}{5.5in}
\end{minipage}
\end{center}

\section{Introduction}\noindent

In this paper, we investigate the decay properties of axially symmetric solutions to the steady Magnetohydrodynamics equations (MHD), which are listed as follows.
\begin{equation}\lab{s-mhd}
\left\{
             \begin{array}{l}
             \textbf{u}\cdot\nabla \textbf{u}+\nabla p=\textbf{h}\cdot\nabla \textbf{h}+\triangle \textbf{u}+\textbf{f},  \\
             \textbf{u}\cdot\nabla \textbf{h}-\textbf{h}\cdot\nabla \textbf{u}=\triangle \textbf{h}+\rm{curl} \,\textbf{g},\\
             \na\cdot\textbf{u} = \na \cdot \textbf{h}=0,
             \end{array}
\right.
\end{equation}
with the additional condition at infinity
\be\lab{infinity condition}
             \lim\limits_{|x|\to\infty}\textbf{u}(x)=\lim\limits_{|x|\to\infty}\textbf{h}(x)=0,
\ee
and the finite Dirichlet integral
\begin{equation}\lab{dirichlet integral}
\int_{\mathbb{R}^3}|\nabla\textbf{u}(x)|^2 + |\na {\bf h} (x)|^2 dx<+\infty.
\end{equation}

Here $\bf{u}$ is the velocity field, $\bf{h}$ is the magnetic field, $p$ is the pressure, $\bf{f}$ and $\na \times {\bf g}$ are the external forces on the magnetically charged fluid flows. Set $\bar{\g} \= \curl \g$ and $\bar{\ff} \= \curl \ff $ in the following. It is well known that Leray \cite{leray33} constructed the weak solutions to the stationary Navier-Stokes equations with no-slip boundary conditions and constant velocity at infinity. Leray's solution has finite Dirichlet integral, and is usually refereed as D-solution. One may refer to the reference \cite{farwig98,finn59,finn65,fujita,lady69} for different construction methods and solutions in different kinds of function spaces. One can show that any weak solution to (\ref{s-mhd}) satisfying (\ref{infinity condition}) and (\ref{dirichlet integral}) is smooth, by following the argument developed in \cite{galdi11}.

Recently there are several papers working on decay rate estimates for the axisymmetric D-solutions to the steady Navier-Stokes equations. Our goal here is to extend the results in \cite{CPZ18,cj9,weng15} to the axisymmetric steady MHD equations. Note that the weighted energy method developed in \cite{cj9,weng15} needs to use the special structure of the vorticity equations, such an extension may be quite nontrivial, although the steady MHD equations has the same scaling as the steady Navier-Stokes equations, the scaling technique developed in \cite{CPZ18} still works in the MHD setting.

Let us first introduce the cylindrical coordinate
\be\no
r=\sqrt{x^2_1+x_2^2},\q \tan \theta=\frac{x_2}{x_1},\q z=x_3,
 \ee
and the basis vectors $\textbf {e}_r$, $\textbf {e}_{\theta}$, $\textbf {e}_z$ are
\be\no
{\textbf {e}_r}=(\cos\,\theta, \sin\,\theta,0),\quad {\textbf {e}_{\theta}}=(-\sin\,\theta,\cos \,\theta ,0),\quad {\textbf {e}_z}=(0,0,1).
\ee
 For convenience, we write $\textbf{x}^{\prime}=(x_1,x_2)$ and $x_3=z.$

We call a function $f$ is \textit{axially symmetric} if it does not depend on $\th$. A vector-valued function $\textbf{u} = (u_r, u_{\th}, u_z)$ is called \textit{axially symmetric} if $u_r,\, u_{\th}, \, u_z$ do not depend on $\th$. And we call a vector-valued function $\textbf{u} = (u_r, u_{\th}, u_z)$ is \textit{axially symmetric without swirl} if $u_{\th}=0$ while $u_r$ and $u_z$ do not depend on $\th$.

Assume that the vector field
$\textbf{u}({\bf x})=u_r(r,z){\textbf {e}_r}+u_{\theta}(r,z){\textbf {e}_{\theta}}+u_z(r,z){\textbf {e}_z}$ and the magnetic field $\textbf{h}({\bf x})=h_r(r,z){\textbf {e}_r}+h_{\theta}(r,z){\textbf {e}_{\theta}}+h_z(r,z){\textbf {e}_z}$ is the axially symmetric solution of (\ref{s-mhd}), then we have the following asymmetric steady MHD equations,
\begin{equation}\lab{asymmetric steady mhd}
\left\{
             \begin{array}{l}
             (u_r\partial_r+u_z\partial_z)u_r-\frac{u^2_{\theta}}{r}+\partial_rp
             =(h_r\partial_r+h_z\partial_z)h_r-\frac{h^2_{\theta}}{r}+(\partial^2_r+\frac{1}{r}\partial_r+\partial^2_z-\frac{1}{r^2})u_r+f_r,   \\
             (u_r\partial_r+u_z\partial_z)u_{\theta}+\frac{u_ru_{\theta}}{r}
             =(h_r\partial_r+h_z\partial_z)h_{\theta}+\frac{h_rh_{\theta}}{r}+(\partial^2_r+\frac{1}{r}\partial_r+\partial^2_z-\frac{1}{r^2})u_{\theta}+f_{\theta},\\
             (u_r\partial_r+u_z\partial_z)u_z+\partial_zp=(h_r\partial_r+h_z\partial_z)h_z+(\partial^2_r+\frac{1}{r}\partial_r+\partial^2_z)u_z+f_z, \\
             (u_r\partial_r+u_z\partial_z)h_r
             =(h_r\partial_r+h_z\partial_z)u_r+(\partial^2_r+\frac{1}{r}\partial_r+\partial^2_z-\frac{1}{r^2})h_r-\partial_zg_{\theta},\\
             (u_r\partial_r+u_z\partial_z)h_{\theta}+\frac{h_ru_{\theta}}{r}=(h_r\partial_r+h_z\partial_z)u_{\theta}+\frac{u_rh_{\theta}}{r}
             +(\partial^2_r+\frac{1}{r}\partial_r+\partial^2_z-\frac{1}{r^2})h_{\theta}+(\partial_zg_r-\partial_rg_z),\\
             (u_r\partial_r+u_z\partial_z)h_z
             =(h_r\partial_r+h_z\partial_z)u_z+(\partial^2_r+\frac{1}{r}\partial_r+\partial^2_z)h_z+\frac{1}{r}\partial_r(rg_{\theta}),\\
             \partial_ru_r+\frac{u_r}{r}+\partial_zu_z=0,\\
             \partial_rh_r+\frac{h_r}{r}+\partial_zh_z=0.
             \end{array}
\right.
\end{equation}

The vorticity $ {\bf \omega^u}$ and the current density ${\bf \omega^h}$ are defined as ${\bf \omega^u}({\bf x})=\nabla\times {\bf u(x)}=\omega_r^u(r,z){\textbf {e}_r}+\omega_{\theta}^u(r,z){\textbf {e}_{\theta}}+\omega_z^u(r,z){\textbf {e}_z}$ and ${\bf \omega^h}(x)=\nabla\times {\bf h(x)}=\omega_r^h(r,z){\textbf {e}_r}+\omega_{\theta}^h(r,z){\bf e_{\theta}}+\omega_z^h(r,z){\textbf {e}_z}$, where
\be\no\lab{general-vorticity}
&\omega_r^u(r,z)=-\partial_zu_{\theta},\;\omega_{\theta}^u=\partial_zu_r-\partial_ru_z,\; \omega_z^u=\frac{1}{r}\partial_r(ru_{\theta}),\\
&\omega_r^h(r,z)=-\partial_zh_{\theta},\;\om_{\th}^h=\p_z h_r -\p_r h_z,\; \omega_z^h=\frac{1}{r}\partial_r(rh_{\theta}).
 \ee
Then we have
\begin{eqnarray}\lab{vorticity eq}
\left\{
\begin{array}{l}
(u_r\partial_r+u_z\partial_z)\omega_r^u-(\omega_r^u\partial_r+\omega_z^u\partial_z)u_r
=(h_r\partial_r+h_z\partial_z)\omega_r^h-(\omega_r^h\partial_r+\omega_z^h\partial_z)h_r\\
\quad\quad\quad+(\partial^2_r+\frac{1}{r}\partial_r+\partial^2_z-\frac{1}{r^2})\omega_r^u- \p_z f_{\th},\\
(u_r\partial_r+u_z\partial_z)\omega_{\theta}^u-\frac{u_r\omega_{\theta}^u}{r}
+2\frac{u_{\theta}}{r}\omega_r^u
=(h_r\partial_r+h_z\partial_z)\omega_{\theta}^h-\frac{h_r\omega_{\theta}^h}{r}
+2\frac{h_{\theta}}{r}\omega_r^h\\
\quad\quad\quad+(\partial^2_r+\frac{1}{r}\partial_r+\partial^2_z-\frac{1}{r^2})\omega_{\theta}^u+ ( \p_z f_r - \p_r f_z ),\\
(u_r\partial_r+u_z\partial_z)\omega_z^u-(\omega_r^u\partial_r+\omega_z^u\partial_z)u_z
=(h_r\partial_r+h_z\partial_z)\omega_z^h-(\omega_r^h\partial_r+\omega_z^h\partial_z)h_z\\
\quad\quad\quad+(\partial^2_r+\frac{1}{r}\partial_r+\partial^2_z)\omega_z^u + \f 1r \p_r (r f_{\th}),\\
(u_r\partial_r+u_z\partial_z)\omega_r^h
-(\omega_r^h\partial_r+\omega_z^h\partial_z)u_r
+\frac{2}{r}\partial_z(u_rh_{\theta})=(h_r\partial_r+h_z\partial_z)\omega_r^u-(\omega_r^u\partial_r+\omega_z^u\partial_z)h_r\\
\quad\quad\quad
+\frac{2}{r}\partial_z(u_{\theta}h_r)+(\partial^2_r+\frac{1}{r}\partial_r+\partial^2_z-\frac{1}{r^2})\omega_r^h - \p_z ( \p_z g_r - \p_r g_z ),
\\
(u_r\partial_r+u_z\partial_z)\omega_{\theta}^h
-\frac{u_r}{r}\omega^h_{\theta}
+2\partial_zu_r\partial_rh_r+2\partial_zu_z\partial_rh_z
=(h_r\partial_r+h_z\partial_z)\omega_{\theta}^u
-\frac{h_r}{r}\omega_{\theta}^u\\
\quad\quad\quad
+2\partial_zh_r\partial_ru_r+2\partial_zh_z\partial_ru_z+(\partial^2_r+\frac{1}{r}\partial_r+\partial^2_z-\frac{1}{r^2})\omega_{\theta}^h - \p_z (\p_z g_{\th}) - \p_r (\f 1r \p_r (r g_{\th})) ,
\\
(u_r\partial_r+u_z\partial_z)\omega_z^h
+(\p_r u_r \p_r + \p_r u_z \p_z)h_{\th}
+ \p_r(\f {u_{\th} h_r}{r})
= (h_r\partial_r+h_z\p_z)\omega_z^u\\ \quad\quad\quad
+(\p_r h_r \p_r + \p_r h_z \p_z)u_{\th}
+ \p_r(\f {u_r h_{\th}}{r})
+(\partial^2_r+\frac{1}{r}\partial_r+\partial^2_z)\omega_z^h + \f 1r \p_r ( r (\p_z g_r - \p_r g_z)).
\end{array}
\right.
\end{eqnarray}

Finn \cite{finn65} have studied the physically reasonable solutions to the steady Navier-Stokes with $O(|x|^{-1})$ decay at infinity. Borchers-Miyakawa \cite{bm95}, Galdi-Simader \cite{gs94} and Novotny-Padula \cite{np95} obtained the existence of such ${\bf u}$ with decay rate $O(|x|^{-1})$ at infinity provided the forcing term is sufficiently small in different function spaces. Under the assumption that $|{\bf u}(x)|=O(|x|^{-1})$ as $|x|\to \oo$, Sverak and Tsai \cite{st2000} had showed that $|\na^k {\bf u}(x)|\leq O(|x|^{-k-1})$ for any $k\geq 1$ by employing a scaling argument and an interior estimate for the Stokes system. In \cite{dg00}, Deuring and Galdi showed that when these solutions are asymptotically expanded near infinity, the leading term cannot be the product of a non-zero vector with the Stokes fundamental solution. One may refer to \cite{np00} for a decomposition of the asymptotic profile at infinity. Korolev and Sverak \cite{ks11} showed that Landau's solution is the right leading term describing the asymptotic behavior when the solution satisfied some smallness assumptions.

The decay properties of the smooth D-solutions to the steady Navier-Stokes and MHD equations are closely related to the uniqueness and energy conservation of the D-solutions. There is a famous open problem which concerns whether ${\bf u}\equiv 0$ is a unique D-solution to steady Navier-Stokes equations. Galdi \cite{galdi11} showed that the D-solution must be trivial if ${\bf u} \in L^{9/2}(\mbR^3)$, which was logarithmically improved in \cite{cwo16}. Many authors have identified different integrability or decay conditions on ${\bf u}$, which lead to many interesting Liouville type results, one may refer to \cite{cpzz19,chae14,cwe16,cw19,chae20,kpr15,ktw17,seregin16} for more details. The Liouville theorem in the case of the axisymmetric three dimensional flows without swirl can be derived from \cite{KNSS09}. Chae \cite{chae14} explored the maximum principle for the total head pressure, and proved the triviality of ${\bf u}$ by assuming $ \De {\bf u} \in L^{\f 65} (\mbR^3)$. Seregin \cite{seregin16,seregin19} applied Caccioppoli type inequality to show that any smooth solution to the stationary Navier-Stokes system in $\mbR^3$, belonging to $L^6(\mbR^3)$ and $BMO^{-1}$, must be zero. Kozono, Terasawa and Wakasugi \cite{ktw17} proposed a decay condition on the vorticity to guarantee that ${\bf u}\equiv 0$. For the MHD system, Liu and Zhang \cite{lz18} proved the Liouville theorem for the bounded ancient solution. In \cite{ww19}, it was shown that there are not non-trivial solutions to MHD equations under the finite Dirichlet integral or the $L^p$ integrability of the velocity and magnetic fields. In their proof, the smallness of the magnetic field plays a vital role.

For the investigation of the decay properties of D-solutions to the exterior stationary Navier-Stokes equations with large external force, Gilbarg and Weinberger \cite{gw74} had made great progress in the two dimensional exterior domain, showed that the weak solution constructed by Leray \cite{leray33} was bounded and converged to a limit $\textbf{u}_0$ in a mean square sense, while the pressure converged pointwise. In \cite{gw78}, they further found that the weak solution with finite Dirichlet integral may not be bounded, but it must grow more slowly than $(\ln r)^{1/2}$. By further adapting the ideas in \cite{gw74,gw78} to the 3D axisymmetric setting, the authors in \cite{cj9,weng15} obtained some decay rates for smooth axially symmetric solutions to steady Navier-Stokes equations. On the other hand, it is well-known that if $(\textbf{u},p)$ solves stationary Navier-Stokes equations, so does $(\textbf{u}_{\la}, p_{\la})$ for all $\la > 0$, where $\textbf{u}_{\la}(\textbf{x}) = \la \textbf{u} (\la \textbf{x})$ and $p_{\la} (\textbf{x}) = \la^2 p(\la \textbf{x})$. According to this, the authors in \cite{CPZ18} utilized the \textit{Brezis-Gallouet inequality} to improve the decay rate of the vorticity in \cite{cj9,weng15}. Recently, the decay properties of the second order derivatives were investigated in \cite{lw20}. Combining the results in \cite{CPZ18,cj9,weng15}, one has: for $\textbf{x} \in \mathbb{R}^3$,
\be\no
&&|\textbf{u}(\textbf{x})| \le C (\f {\ln r}{r})^{1/2},\ \ \ \ |\om_r^u(\textbf{x})|+|\om^u_z (\textbf{x})| \le C \f {(\ln r)^{11/8}}{r^{9/8}},\\\no
&&|\om_{\th}^u(\textbf{x})|+ |\na \om_{\th}^u(\textbf{x})|  \le C \f {(\ln r)^{3/4}}{r^{5/4}},\\\no
&&|\na u_r| + |\na u_z| \le C \f {(\ln r)^{7/4}}{r^{5/4}},
\ee
where $C$ is a positive constant.

As was noticed in \cite{CPZ18}, the derivation of the decay of the velocity itself almost had no use of the Navier-Stokes equations, it is straightforward to gain the same decay rates for the magnetic field. However, the estimates of the vorticity does depend on the structure of the vorticity equations. For the steady MHD equations, the equations for the electric current $\nabla \times {\bf h}$ are more complicate than those of the vorticity, therefore the decay rates of ${\bf \om^u}$ and ${\bf \om^h}$ obtained here are weaker than those of Navier-Stokes equations. Our first result is stated as follows.
\bt\label{main11}
{\it
Support $\ff \in H^1 (\mbR^3)$ and $\g \in H^2 (\mbR^3)$ be axially symmetric vector field with
 \be\lab{fg}
 \|(\ff,\curl \g)\|_{L^{\frac{3}{2}}(\mbR^3)}+\|\textbf{f} \|_{H^1 (\mbR^3)} + \| r^{\f 12} (\ff, \curl \g)\|_{L^2 (\mbR^3)} + \| \textbf{g} \|_{H^2 (\mbR^3)}  \leq M
 \ee
for some constant $M$. Let $(\textbf{u},\textbf{h},p)$ be a smooth axially symmetric solution to the steady MHD equations (\ref{s-mhd}) with (\ref{infinity condition})-(\ref{dirichlet integral}), then there exist a constant $C (M) >0$, such that
\be\lab{uhdecay}
&&|\textbf{u}| + |\textbf{h}|\le C(M)(\f {\ln r}{r})^{\f 12},\\\lab{omthetadecay}
&&|\omega_{\theta}^u| + |\omega_{\theta}^h| \le C(M) \frac{\ln r}{r},\\\lab{omrzdecay}
&&|\omega_r^u| + |\omega_z^u| + |\omega_r^h|+ |\omega_z^h| \le C(M) \frac{\ln r}{r},
\ee
for large $r$.
}\et

Next we consider a special class of axisymmetric D-solutions to the steady MHD equations with $\textbf{h}(r,z) = h_{\th}(r,z) \boldsymbol{e_{\theta}}$. In this situation, we obtain better decay results for the vorticity field by employing the \textit{Biot-Savart law} as in \cite{CPZ18}. However, this argument does not improve the decay rates of $\om_r^h$ and $\om_z^h$. Here we resort to the weighted energy estimates in \cite{weng15} to get a better decay for $\om_r^h$ and $\om_z^h$.
\bt\label{main12}
{\it
Support $\ff \in H^2 (\mbR^3)$ and $\g \in H^3 (\mbR^3)$ be axially symmetric vector field without swirl, and satisfying (\ref{fg}) and
 \be\lab{fg1}
 && \| (\f r {\ln r})^{\f 12} f_{\th} \|_{L^2 (\mbR^3)} + \| r (f_r, f_z) \|_{L^2 (\mbR^3)} + \| r^2 \curl \emph{\textbf{f}} \|_{L^2 (\mbR^3)} \le M,\\\lab{fg2}
 &&   \| r^2 \p_r (\p_z g_r - \p_r g_z) \|_{L^2 (\mbR^3)} \le M,
 \ee
 for some constant $M$.
Let $(\textbf{u},\textbf{h},p)$ be a smooth axially symmetric solution to the steady MHD equations (\ref{s-mhd}) with (\ref{infinity condition})-(\ref{dirichlet integral}), where $\textbf{h}(r,z)=h_{\theta}(r,z)\boldsymbol{e_{\theta}}$. Then we have better estimates as follows:
\be\lab{decay-om-th-u}
&&|\omega_{\theta}^u|  \leq C(M) \frac{(\ln r)^{3/4}}{r^{5/4}},\\\lab{decay-om-rz-u}
&&|\omega_r^u|+|\omega_z^u| \leq C(M) \frac{(\ln r)^{11/8}}{r^{9/8}},\\\lab{decay-na-u-rz}
&&|\nabla u_r|+|\nabla u_z| \leq  C(M) \frac{(\ln r )^{7/4}}{r^{5/4}},\\\lab{decay-na-h-th}
&&|\om^h_r|+ |\om^h_z| + |\na h_{\th} | \leq C(M) r^{- (\f {35}{32})^-},
\ee
where $a^-$ denotes any number less than $a$.
}\et

The last result examines the decay rate in the $Oz$-direction for the axisymmetric D-solution of steady MHD equations with the magnetic field having only the swirl component ${\bf h} (r,z) = h_{\th} (r,z) {\bf e}_{\th}$. The quantity $\Pi \=\frac{h_{\theta}}{r}$ satisfies an elliptic equation, from which weighted energy estimates are available for the weight $\rho=\sqrt{r^2+z^2}$. Furthermore, if the swirl of the velocity $u_{\th}\equiv 0$, we can deduce some decay estimates for the velocity as in steady axisymmetric Navier-Stokes case \cite{weng15}. However, the decay rates are worse than the Navier-Stokes case.
\bt\label{main13}
{\it
 Suppose $ \ff \in H^1 (\mbR^3)$ and $ \g \in H^2 (\mbR^3)$ be axially symmetric vector field without swirl, and satisfying (\ref{fg}), (\ref{fg2}) and
\be\lab{fg3}
&&\| \rho \f {\p_z f_r - \p_r f_z}{r} \|_{L^2 (\mbR^3)} \le M,\\
\lab{fg4}
&&\| \rho \f {\p_z g_r - \p_r g_z}{r} \|_{L^2 (\mbR^3)} \le M\\\lab{fg5}
&& \| (\f r {\ln r})^{\f 12} g_{\th} \|_{L^2 (\mbR^3)} + \| r^2 \curl \emph{\textbf{g}} \|_{L^2 (\mbR^2)} \le M,
\ee
for some constant $M$.
Let $(\textbf{u},\textbf{h},p)$ be a smooth axially symmetric solution to the steady MHD equations (\ref{s-mhd}) with (\ref{infinity condition})-(\ref{dirichlet integral}), where $\textbf{h}(r,z)=h_{\theta}(r,z)\boldsymbol{e_{\theta}}$.
Then there holds
\be\no
|h_{\th} (r,z)| \le C(M) (\rho + 1)^{- \f {31}{168}}, \quad \quad  \rho=\sqrt{r^2+z^2}.
\ee

Moreover, if the swirl velocity $u_{\th}\equiv 0$, then
\be
|u_r(r,z)| +|u_z(r,z)|&\leq& C(M) (\rho+1)^{-(\f {13}{179})^-},\\\lab{omgea300}
|\om_{\th}(r,z)| &\leq& C(M)(\rho+1)^{-(\f {39}{179})^-},\\\lab{htheta30}
|h_{\th}(r,z)|&\leq & C(M) (\rho + 1)^{-(\f {3041}{15036})^-}.
\ee
}\et

This paper is organized as follows. In section 2 we introduce some preliminary tools. In Section 3, we prove the decay rates of $({\bf u}, {\bf h})$ and $({\bf \om^u}, {\bf \om^h})$ by applying the scaling argument and the \textit{Brezis-Gallouet inequality}. In Section 4, we will use the \textit{Biot-Savart law} and also weighted energy estimates for the magnetic field to get some improved decay rate estimates for the special class of solutions with ${\bf h}(r,z) = h_{\th} (r,z) {\bf e}_{\th}$. In Section 5, we will investigate the decay properties of $({\bf u}, {\bf h})$ in the $O_z$-direction also within the special class of solutions with ${\bf h}(r,z) = h_{\th} (r,z) {\bf e}_{\th}$.

\section{Preliminary}

\bl\lab{BG ineq}\emph{\cite[Theorem 1]{bj80}}
{\it Let $f\in H^2(\Omega)$ where $\Omega\subset\mathbb{R}^2$. Then there exists a constant $C_{\Omega}$, depending only on $\Omega$, such that
\begin{equation}\no
\| f \|_{L^{\infty}(\Omega)}\le C_{\Omega}(1+\| f \|_{H^1(\Omega)})\ln^{1/2}\big(e+\| \Delta f \|_{L^2(\Omega)}\big).
\end{equation}
}
\el

\bl\lab{lemma CZ}\emph{\cite[Lemma 3.2]{CPZ18}}
{\it Assume that $K(x,y)$ be a Calderon-Zygmund kernel and $f$ is a smooth axisymmetric function satisfying, for ${\bf x}=({\bf x^{\prime}}, z) \in \mbR^3$
\be\no
|f(x)| + |\na f(x)| \le \f {\ln^b (e+ |x^{\prime}|)}{(1+ |x^{\prime}|)^a} \q for \q 0<a<2,\q b>0.
\ee
Define $Tf(x) \= \int K(x,y)f(y) dy$. Then there exists a constant $c_0$ such that
\be\no
|Tf(x)|\le c_0 \f {\ln^{b+1} (e+ |x^{\prime}|)}{(1+ |x^{\prime}|)^a}.
\ee
}\el

\bl\lab{a decay lemma}\emph{\cite[Lemma 2.1]{weng15}}
{\it Suppose a smooth axially symmetric function $f(x)$ satisfies the following weighted energy estimates
\begin{equation}\nonumber
\int_{\mathbb{R}^3}\left( r^{e_1}|f(r,z)|^2+ r^{e_2}|\nabla f(r,z)|^2 + r^{e_3}|\partial_z\nabla f(r,z)|^2\right)dx\leq C,
\end{equation}
with nonnegative constants $e_1$, $e_2$, $e_3$. Then for any $r>0$, we have
\begin{eqnarray}\no
&&\int_{-\infty}^{\infty}|f(r,z)|^2dz\leq Cr^{-\frac{1}{2}(e_1+e_2)-1},\\\no
&&\int_{-\infty}^{\infty}|\partial_zf(r,z)|^2dz\leq Cr^{-\frac{1}{2}(e_3+e_2)-1},\\\no
&&|f(r,z)|^2\le Cr^{-\frac{1}{4}(e_1+2e_2+e_3)-1}, \quad \text{$\forall z\in\mathbb{R}^3$}.
\end{eqnarray}
}\el

\section{Proof of Theorem \ref{main11}} \noindent

Before starting to prove Theorem \ref{main11}, we want to sketch the proof as follows. Following \cite{CPZ18}, the decay rate of $(\u, \h)$ can be shown by a scaling argument and the \textit{Brezis-Gallouet} inequality. Next utilizing the localized energy estimates of $\om^u$ and $\om^h$ and the $\textit{Brezis-Gallout}$ inequality, it will be deduced that
\be\no
& &
{|\om_{\th}^u (r,z) | + |\om_{\th}^h (r,z) | \le}
{C(M ) (\f {\ln r}r)^{\f 12} \| (\u, \h) \|_{L^{\oo} ([\f r2, \f {3r}2] \times [-r, r])},}\\\no
& &
{|\om_r^u (r,z)| + |\om_z^u (r,z)| + |\om_r^h (r,z) | + |\om_z^h (r,z)|}
{\le}
{C (M) (\f {\ln r}r)^{\f 12} \| (\u, \h) \|_{L^{\oo} ([\f r2, \f {3r}2] \times [-r, r])}.}
\ee
The details can be found in (\ref{lam-rz})  and (\ref{lam-th}). The estimates \eqref{omthetadecay} and \eqref{omrzdecay} follow directly from the above two inequalities.

Now we start to prove Theorem \ref{main11}. Assume that \eqref{fg} holds, the standard existence theory tells us that there exists a weak solution $({\bf u}, {\bf h} ,p)$ to MHD equation (1.1)-(1.2) with the finite Dirichlet integral (1.3) (see Chapter X in \cite{galdi11} for more details). Since $({\bf u}, {\bf h})\in L^6(\mbR^3)$ and $(\nabla {\bf u},\nabla {\bf h})\in L^2(\mbR^3)$, then $({\bf u}\cdot\nabla {\bf u},{\bf h}\cdot\nabla {\bf h}, {\bf u}\cdot\nabla {\bf h},{\bf h}\cdot\nabla {\bf u})\in L^\frac{3}{2}(\mbR^3)$.  According to the $L^p$ estimates to the Stokes system (See Chapter IV.2 in \cite{galdi11}), if $(\ff,\text{curl }{\bf g})\in  L^{\frac{3}{2}} (\mbR^3)$, then $(\nabla^2 {\bf u}, \nabla^2{\bf h})\in L^{\frac{3}{2}} (\mbR^3)$, which implies that $(\nabla {\bf u},\nabla {\bf h})\in L^3(\mathbb{R}^3)$ by Sobolev embedding. Thus $({\bf u}\cdot\nabla {\bf u},{\bf h}\cdot\nabla {\bf h}, {\bf u}\cdot\nabla {\bf h},{\bf h}\cdot\nabla {\bf u})\in L^2(\mbR^3)$, yielding that
\be\label{dd}
\|(\na^2 \u,\na^2 \h)\|_{L^2(\mathbb{R}^3)} \le C(M),
 \ee
where the positive constant $C(M)$ depending only on $M$. And there also hold $\|({\bf u}, {\bf h})\|_{W^{1,6}(\mathbb{R}^3)}+\|({\bf u}, {\bf h})\|_{L^{\oo} (\mbR^3)} \leq C (M)$. Simple calculations show that $({\bf u}\cdot\nabla {\bf u},{\bf h}\cdot\nabla {\bf h}, {\bf u}\cdot\nabla {\bf h},{\bf h}\cdot\nabla {\bf u})\in H^1(\mbR^3)$. Then the $L^p$ theory to the Stokes system yields that
\be\label{third}
\|(\na^3 \u,\na^3 \h)\|_{L^2(\mathbb{R}^3)} \leq C(M).
\ee

Fix any point ${\bf x_0} \in\mathbb{R}^3$ such that $|{\bf x_0^{\prime}}|=\lambda$ is large, and consider the scaled solutions $\tilde{\textbf{u}}(\tilde{\textbf{x}})$, $\tilde{\textbf{h}}(\tilde{\textbf{x}})$ and the two dimensional domain $\tilde{D}$,
\be\no
&&\tilde{\textbf{u}}(\tilde{\textbf{x}})=  \la {\textbf{u}}(\la\tilde{\textbf{x}}),\ \ \ \tilde{\textbf{h}}(\tilde{\textbf{x}})= \la {\textbf{h}}(\la\tilde{\textbf{x}}),\\\no
&&\tilde{D}=  \{(\tilde{r}, \tilde{z})|1/2\le \tilde{r}\le2,|\tilde{z}|\le1\},
\ee
where ${\textbf{x}}= \la \tilde{\textbf{x}}$.

To prove the estimate \eqref{uhdecay}, by carefully checking the argument developed in \cite{CPZ18} for the steady axisymmetric Navier-Stokes equations without external force, only the estimates \eqref{dd} and the finite Dirichlet integral are used in their proof, the structure of equations play no role, thus the estimate \eqref{uhdecay} holds.

It remains to derive the decay rates of $\bf{\omega^u}$ and $\bf{\omega^h}$. According to previous scaling, we have
\be\no
&\tilde{\om}^u (\tilde{\textbf{x}})= \lambda^2 \om^u(\lambda \tilde{\textbf{x}}) = \lambda^2 \om^u({\textbf{x}}),\\\no
&\tilde{\om}^h (\tilde{\textbf{x}})= \lambda^2 \om^h(\lambda \tilde{\textbf{x}}) = \lambda^2 \om^h({\textbf{x}}),\\\no
& \tilde{\ff} (\tilde{\x}) = \la^3 \ff (\la \tilde{\x}) = \la^3 \ff (\x),\\\no
&  \curl_{\tilde{\x}} \tilde{\g} (\tilde{\x}) = \la^3 \curl_{\x} \g (\la \tilde{\x}) = \la^3 \curl_{\x} \g (\x).
\ee
For simplification of notation, we will drop the $`` \sim "$ when computations take place under the scaled sense. Select the domains
\be\no
&\mathcal{C}_1 = \{ (r, \th, z) : \f 12 < r < \f 32, 0 \le \th \le 2\pi, |z| \le 1 \},\\\no
&\mathcal{C}_2 = \{ (r, \th, z) : \f 34 < r < \f 54, 0 \le \th \le 2\pi, |z| \le \f 12 \}.
\ee

Let $\phi(y)$ be a cut-off function satisfying supp $\phi(y) \subset \mathcal{C}_1$ and $\phi(y)=1$ for $y\in\mathcal{C}_2$ such that the gradient of $\phi$ is bounded. Now testing the vorticity equation (\ref{vorticity eq}) with $\omega_r^u\phi^2$, $\omega_{\theta}^u\phi^2$, $\omega_z^u\phi^2$, $\omega_r^h\phi^2$, $\omega_{\theta}^h\phi^2$ and $\omega_z^h\phi^2$ respectively, and integrating over $\mathcal{C}_1$, after direct computations we obtain
\be\no
&&\int_{\mathcal{C}_1}\bigg(|\nabla(\omega_{\theta}^u\phi)|^2+\frac{(\omega_{\theta}^u\phi)^2}{r^2}\bigg)dy=\int_{\mathcal{C}_1}\bigg(|\omega_{\theta}^u|^2|\nabla\phi|^2
+\frac{1}{2}|\omega_{\theta}^u|^2(u_r\partial_r+u_z\partial_z)(\phi^2)\\\no
&&\quad\quad+ \f {u_r}{r} (\om_{\th}^u \phi)^2
-\frac{2}{r}u_{\theta}\omega_r^u\omega_{\theta}^u\phi^2
+\frac{2}{r}h_{\theta}\omega_r^h\omega_{\theta}^u\phi^2+\omega_{\theta}^u\phi(h_r\partial_r+h_z\partial_z)(\omega_{\theta}^h\phi)-\frac{h_r}{r}\omega_{\theta}^u\omega_{\theta}^h\phi^2\\\no
&&\quad\quad-\omega_{\theta}^u\omega_{\theta}^h\phi(h_r\partial_r+h_z\partial_z)\phi
+ \phi(f_z \p_r - f_r \p_z) (\om_{\th}^u \phi) + \om_{\th}^u \phi (f_z \p_r - f_r \p_z) \phi+ f_z \f {\om_{\th}^u \phi}{r} \phi
\bigg)dy\\\no
&&\leq
C
(1+\|({\bf u}, {\bf h})\|^2_{L^{\infty}(\mathcal{C}_1)})\|(\omega_r^u,\omega_{\theta}^u,\omega_r^h,\omega_{\theta}^h)\|^2_{L^2(\mathcal{C}_1)}
+ C \| (f_r, f_z ) \|_{L^2 (\mc{C}_1)}^2
\\\no
&& \q
+\frac{1}{8}\| (\nabla(\omega^u_{\theta}\phi), \nabla(\omega^h_{\theta}\phi))\|^2_{L^2(\mathcal{C}_1)},
\\\no
&&\int_{\mathcal{C}_1}\bigg(|\nabla(\omega_{\theta}^h\phi)|^2+\frac{(\omega_{\theta}^h\phi)^2}{r^2}\bigg)dy=\int_{\mathcal{C}_1}\bigg(
|\na \phi|^2 |\om_{\th}^h|^2 + \f 12 |\om_{\th}^h|^2 (u_r \p_r + u_z \p_z) (\phi^2) \\\no
&&\quad\quad+ 2 \p_z(\om_{\th}^h \phi^2) (u_z \p_r h_z - h_r \p_r u_r)+ 2 \p_r(\om_{\th}^h \phi^2) (h_r \p_z u_r - u_z \p_z h_z)+ \f{u_r}{r} (\om_{\th}^h \phi)^2\\\no
&&\quad\quad+ \om_{\th}^h \phi (h_r \p_r + h_z \p_z) (\om_{\th}^u \phi) - \om_{\th}^u \om_{\th}^h \phi (h_r \p_r + h_z \p_z) \phi- \f {h_r}{r} \om_{\th}^u \phi \cdot \om_{\th}^h \phi \\\no
&&\quad\quad+ 2 \f {\om_{\th}^h \phi}{r} (h_r \p_z u_r - u_z \p_z h_z)\phi + \p_z g_{\th}(\p_z (\om_{\th}^h \phi))+ \om_{\th}^h  \phi \p_z \phi) \\\no
&& \q\q  + \p_r (r g_{\th}) (\phi \p_r (\om_{\th}^h \phi) + \om_{\th}^h \phi \p_r \phi + \f {\om_{\th}^h \phi}{r} \phi )
\bigg)dy\\\no
&&\leq
C (1+\|(u_r,u_z,h_r,h_z)\|^2_{L^{\infty}(\mathcal{C}_1)})\|(\omega_{\theta}^u,\omega_{\theta}^h)\|^2_{L^2(\mathcal{C}_1)}
+ C \| \bar{\g} \|_{L^2 (\mc{C}_1)}^2
\\\no
&&\quad\quad+ C\|(u_r,u_z,h_r,h_z)\|^2_{L^{\infty}(\mathcal{C}_1)}\|(\nabla u_r,\nabla h_z)\|^2_{L^2(\mathcal{C}_1)}+\frac{1}{8}\|(\na (\om_{\th}^u \phi),\nabla(\omega_{\theta}^h\phi))\|^2_{L^2(\mathcal{C}_1)}.
\ee

Combining these two estimates, we find
\be\nonumber
&&\|(\nabla\omega_{\theta}^u,\nabla\omega_{\theta}^h)\|^2_{L^2(\mathcal{C}_2)}\leq C (1+\|(\textbf{u},\textbf{h})\|^2_{L^{\infty}(\mathcal{C}_1)})\|(\omega_r^u,\omega_{\theta}^u,\omega_r^h,\omega_{\theta}^h)\|^2_{L^2(\mathcal{C}_1)}
  \\\lab{C1-om-th}
&&\quad \q
+C \| (f_r, f_z , \bar{\g}) \|_{L^2 (\mc{C}_1)}^2
+
C  \|(u_r,u_z,h_r,h_z)\|^2_{L^{\infty}(\mathcal{C}_1)}\|(\nabla u_r,\nabla h_z)\|^2_{L^2(\mathcal{C}_1)}.
\ee
There also holds
\be\no
&&\int_{\mathcal{C}_1}\bigg(|\nabla(\omega_r^u\phi)|^2+\frac{(\omega_r^u\phi)^2}{r^2}\bigg)dy= \int_{\mathcal{C}_1} \bigg(|\na \phi|^2 |\om_r^u|^2 + \f 12 |\om_r^u|^2 (u_r \p_r + u_z \p_z) (\phi^2) \\\no
&&\quad\quad- \f {u_r}{r} (\om_r^u \phi)^2 - u_r (\p_r (\om_r^u \phi)^2
+\p_z (\om_r^u \phi \cdot \om_z^u \phi)) + \om_r^u \phi (h_r \p_r + h_z \p_z)(\om_r^h \phi)
\\\no &&\quad\quad- \om_r^u \om_r^h \phi (h_r \p_r + h_z \p_z) \phi
+ h_r \p_r(\om_r^u \phi \cdot \om_r^h \phi) + h_r\p_z (\om_r^u \phi \cdot \om_z^h \phi)\\\no
&& \q\q + \f {h_r}{r} \om_r^u\phi \cdot \om_r^h \phi + f_{\th}\phi (\p_z (\om_r^u \phi) + \om_r^u \p_z \phi)
\bigg)dy\\\no
&&\leq C
(1+\|(u_r,u_z,h_r,h_z)\|^2_{L^{\infty}(\mathcal{C}_1)})\|(\omega_r^u,\omega_z^u,\omega_r^h,\omega_z^h)\|^2_{L^2(\mathcal{C}_1)}
+ C \| f_{\th } \|_{L^2 (\mc{C}_1)}^2 \\\no
&&\quad\quad
+\frac{1}{8}\|(\nabla(\omega_r^u\phi),\nabla(\omega_z^u\phi),\nabla(\omega_r^h\phi),\nabla(\omega_z^h\phi)\|^2_{L^2(\mathcal{C}_1)},
\ee
\be\no
&&\int_{\mathcal{C}_1}|\nabla(\omega_z^u\phi)|^2dy= \int_{\mathcal{C}_1}\bigg(
|\na \phi|^2 |\om_z^u|^2 + \f 12 |\om_z^u|^2 (u_r \p_r + u_z \p_z) (\phi^2)
-u_z \p_r (\om_r^u \phi \cdot \om_z^u \phi)\\\no
&&\quad\quad
- \f {u_z}{r} \om_r^u \phi \cdot \om_z^u \phi + \om_z^u \phi (h_r \p_r + h_z \p_z) (\om_z^h \phi)
-\om_z^u \om_z^h \phi (h_r \p_r + h_z \p_z )\phi\\\no
&&\quad\quad
+ h_z (\p_r (\om_z^u \phi \cdot \om_r^h \phi) + \p_z (\om_z^u \phi \cdot \om_z^h \phi))
 - u_z \p_z (\om_z^u \phi)^2+ \f {h_z}{r} \om_z^u \phi \cdot \om_r^h \phi
\\\no
  && \q\q
  - f_{\th} \phi (\p_r (\om_z^u \phi) + \om_z^u \p_r \phi)
\bigg) dy\\\no
&&\leq C
(1+\|(u_r,u_z,h_r,h_z)\|^2_{L^{\infty}(\mathcal{C}_1)})\|(\omega_r^u,\omega_z^u,\omega_r^h,\omega_z^h)\|^2_{L^2(\mathcal{C}_1)}
+ C \| f_{\th } \|_{L^2 (\mc{C}_1)}^2 \\\no
&&\quad\quad
+ \f 18 \|(\nabla(\omega_r^u\phi),\nabla(\omega_z^u\phi),\nabla(\omega_r^h\phi),\nabla(\omega_z^h\phi))\|^2_{L^2(\mathcal{C}_1)},
\\\no
&&\int_{\mathcal{C}_1}\bigg(|\nabla(\omega_r^h\phi)|^2+\frac{(\omega_r^h\phi)^2}{r^2}\bigg)dy=\int_{\mathcal{C}_1}\bigg(
|\na \phi|^2 |\om_r^h|^2 + \f 12 |\om_r^h|^2 (u_r \p_r + u_z \p_z) (\phi^2) \\\no
&&\quad\quad- u_r \p_r(\om_r^h \phi)^2  - u_r \p_z (\om_r^h \phi \cdot \om_z^h \phi)-\f {u_r}{r} |\om_r^h \phi|^2  + \om_r^h \phi (h_r \p_r + h_z \p_z) (\om_r^u \phi) \\\no
&&\quad\quad- \om_r^u \om_r^h \phi (h_r \p_r + h_z \p_z) \phi + h_r \p_r (\om_r^u \phi \cdot \om_r^h \phi)  + h_r \p_z (\om_z^u \phi \cdot \om_r^h \phi)\\\no
&&\quad\quad + \f {h_r}{r} \om_r^u \phi \cdot \om_r^h \phi - \f 2r \p_z (u_r h_{\th}) \om_r^h \phi^2 + \f 2r \p_z (u_{\th} h_r) \om_r^h \phi^2\\\no
&& \q\q  + \p_z (\om_r^h \phi) (\p_z g_r - \p_r g_z)
\bigg)dy\\\no
&&\leq
C
(1+\|(u_r,u_z,h_r,h_z)\|^2_{L^{\infty}(\mathcal{C}_1)})\|(\omega_r^u,\omega_z^u,\omega_r^h,\omega_z^h)\|^2_{L^2(\mathcal{C}_1)}
+ C \| \bar{\g} \|_{L^2 (\mc{C}_1)}^2 \\\no
&&\quad\quad
+
C
\|(u_r,u_{\theta},h_r,h_{\theta})\|^2_{L^{\infty}(\mathcal{C}_1)}
\|\nabla (u_r,u_{\theta},h_r,h_{\theta})\|^2_{L^2(\mathcal{C}_1)}
\\\no
&& \q\q+\frac{1}{8}\|(\nabla(\omega_r^u\phi),\nabla(\omega_z^u\phi),\nabla(\omega_r^h\phi),\nabla(\omega_z^h\phi))\|^2_{L^2(\mathcal{C}_1)}.
\\\no
&&\int_{\mathcal{C}_1}|\nabla(\omega_z^h\phi)|^2dy=\int_{\mathcal{C}_1}\bigg(
|\omega_z^h|^2|\nabla\phi|^2+\frac{1}{2}(\omega_z^h)^2(u_r\partial_r+u_z\partial_z)(\phi^2)\\\no
&&\quad\quad + \p_r(\om_z^h \phi^2) (h_{\th} \p_r u_r -u_{\th} \p_r h_r)+ \p_z(\om_z^h \phi^2) (h_{\th} \p_r u_z - u_{\th} \p_r h_z )\\\no
&&\quad\quad + \f {\om^h_z \phi^2}r (\p_r (u_r h_{\th}) - \p_r (u_{\th} h_r)) + \om_z^h \phi (h_r \p_r + h_z \p_z) (\om_z^u \phi)\\\no
&& \q\q- \om_z^u \om_z^h \phi (h_r \p_r + h_z \p_z)\phi - \p_r (\om_z^h \phi^2) (\p_z g_r - \p_r g_z)
\bigg)dy\\\no
&&\leq C
(1+\|(u_r,u_z,h_r,h_z)\|^2_{L^{\infty}(\mathcal{C}_1)})\|(\omega_z^u,\omega_z^h)\|^2_{L^2(\mathcal{C}_1)}+\frac{1}{8}\|(\nabla(\omega_z^u\phi),\nabla(\omega_z^h\phi))\|^2_{L^2(\mathcal{C}_1)}\\\no
&&\quad\quad
+C
\|(u_r, u_{\theta},h_r, h_{\theta})\|^2_{L^{\infty}(\mathcal{C}_1)}\|(\na \u, \na \h)\|^2_{L^2(\mathcal{C}_1)}
+
C \| \bar{\g}\|_{L^2 (\mc{C}_1)}^2.
\ee
Then we obtain
 \be\nonumber
&&\|(\nabla\omega_r^u,\nabla\omega_z^u,\nabla\omega_r^h,\nabla\omega_z^h)\|^2_{L^2(\mathcal{C}_2)}\leq C  (1+\|(u_r, u_z, h_r, h_z)\|^2_{L^{\infty}(\mathcal{C}_1)})\|(\omega_r^u,\omega_z^u,\omega_r^h,\omega_z^h)\|^2_{L^2(\mathcal{C}_1)}\\\lab{C1-om-rz}
&&\quad\quad +C  \|(u_r,u_{\theta},h_r,h_{\theta})\|^2_{L^{\infty}(\mathcal{C}_1)}
\| (\na \u, \na \h)\|^2_{L^2(\mathcal{C}_1)} +C  \| (f_{\th }, \bar{\g}) \|_{L^2 (\mc{C}_1)}^2 .
\ee
Set
\begin{equation*}
\mathcal{\bar{C}}_2 \= \{(r,z):\frac{3}{4}<r<\frac{5}{4},|z|\le1/2\}.
\end{equation*}
Utilizing the \textit{Brezis-Gallouet inequality} in Lemma \ref{BG ineq} and the localized energy estimates (\ref{C1-om-th}) and (\ref{C1-om-rz}), we can conclude
\be\nonumber
&&\|(\omega_r^u,\omega_z^u,\omega_r^h,\omega_z^h)\|_{L^{\infty}(\mathcal{\bar{C}}_2)}\leq C \bigg(1+
(1+\|(u_r, u_z, h_r, h_z)\|_{L^{\infty}(\mathcal{C}_1)})\|(\omega_r^u,\omega_z^u,\omega_r^h,\omega_z^h)\|_{L^2(\mathcal{C}_1)}
\\\no
&&\quad \quad+\|(u_r,u_{\theta},h_r,h_{\theta})\|_{L^{\infty}(\mathcal{C}_1)}
\|(\na \u, \na \h)\|_{L^2(\mathcal{C}_1)}
+
\| (f_{\th }, \bar{\g}) \|_{L^2 (\mc{C}_1)}
\bigg)\\\no
&&\quad\quad\times
\ln^{1/2}\bigg(e+\|\De(\omega_r^u,\omega_z^u,\omega_r^h,\omega_z^h)\|_{L^2(\bar{\mc{C}}_2)}\bigg),\\\nonumber
&&\|(\omega_{\theta}^u,\omega_{\theta}^h)\|_{L^{\infty}(\mathcal{\bar{C}}_2)}\le C  \bigg(1+
(1+\|(\textbf{u},\textbf{h})\|_{L^{\infty}(\mathcal{C}_1)})\|(\omega_r^u,\omega_{\theta}^u,\omega_r^h,\omega_{\theta}^h)\|_{L^2(\mathcal{C}_1)}
 +
 \| (f_r, f_z , \bar {\g}) \|_{L^2 (\mc{C}_1)}\\\no
&&\q\q+\|(u_r,u_z,h_r,h_z)\|_{L^{\infty}(\mathcal{C}_1)}\|(\nabla u_r,\nabla h_z)\|_{L^2(\mathcal{C}_1)}
\bigg)
\ln^{1/2}\bigg(e+\|(\De\omega_{\theta}^u,\De\omega_{\theta}^h)\|_{L^2(\bar{\mc{C}}_2)}\bigg).
\ee

Then scaling back to the domains
\be\no
\mathcal{C}_{1,\lambda} & = & \{(r,\theta,z):\frac{\lambda}{2}<r<\frac{3\lambda}{2},0\le\theta\le2\pi,|z|\le\lambda\},\\\no
\mathcal{C}_{2,\lambda} & = & \{(r,\theta,z):\frac{3\lambda}{4}<r<\frac{5\lambda}{4},0\le\theta\le2\pi,|z|\le\frac{\lambda}{2}\},
\ee
we derive
\be\nonumber
&&\lambda^2\|(\omega_r^u,\omega_z^u,\omega_r^h,\omega_z^h)\|_{L^{\infty}(\mathcal{C}_{2,\lambda})}\leq
C \bigg( 1 +
(1+\lambda\|(\textbf{u},\textbf{h})\|_{L^{\infty}(\mathcal{C}_{1,\lambda})})\lambda^{1/2}\|(\omega_r^u,\omega_z^u,\omega_r^h,\omega_z^h)\|_{L^2(\mathcal{C}_{1,\lambda})}
\\\no
&&\quad  +\lambda^{3/2}\|(u_r,u_{\theta},h_r,h_{\theta})\|_{L^{\infty}(\mathcal{C}_{1,\lambda})}
\|( \nabla \u, \na \h) \|_{L^2(\mathcal{C}_{1,\lambda})}
+
 \la^{3/2} \| ( f_{\th }, \bar {\g}) \|_{L^2 (\mc{C}_{1, \la})}  \bigg) \\\lab{lam-rz}
 && \q  \times  \ln^{1/2}\bigg(e+\lambda^{5/2}\|\De(\omega_r^u,\omega_z^u,\omega_r^h,\omega_z^h)\|_{L^2(\mathcal{C}_{2,\lambda})}\bigg),
\\\no
&&\lambda^2\|(\omega_{\theta}^u,\omega_{\theta}^h)\|_{L^{\infty}(\mathcal{C}_{2,\lambda})}\le
C \bigg(1 + (1+
\lambda\|(\textbf{u},\textbf{h})\|_{L^{\infty}(\mathcal{C}_{1,\lambda})})
\lambda^{1/2}\|(\omega_r^u,\omega_{\th}^u,\omega_r^h,\omega_{\th}^h)\|_{L^2(\mathcal{C}_{1,\lambda})}
\\\no
&& \quad+\lambda^{3/2}\|(u_z, h_r)\|_{L^{\infty}(\mathcal{C}_{1,\lambda})}\|(\nabla u_r,\nabla h_z)\|_{L^2(\mathcal{C}_{1,\lambda})}
+
 \la^{3/2} \| ( f_r, f_z , \bar {\g}) \|_{L^2 (\mc{C}_{1, \la})}
\bigg)\\\lab{lam-th}
&& \q \times  \ln^{1/2}\bigg(e+\lambda^{5/2}\|(\De\omega_{\th}^u,\De\omega_{\th}^h)\|_{L^2(\mathcal{C}_{2,\lambda})}\bigg).
\ee

By the finite Dirichlet integral assumption (\ref{dirichlet integral}), the estimates \eqref{uhdecay} and \eqref{third}, we conclude that
\begin{eqnarray}\no
&&\|(\omega_{\theta}^u,\omega_{\theta}^h)\|_{L^{\infty}(\mathcal{C}_{2,\lambda})}\leq  C (M ) \frac{\ln \lambda}{\lambda},\\\no
&&\|(\omega_r^u,\omega_z^u,\omega_r^h,\omega_z^h)\|_{L^{\infty}(\mathcal{C}_{2,\lambda})}\leq C (M ) \frac{\ln \lambda}{\lambda},
\end{eqnarray}
 where we use the boundedness of $\| ( \na^3 \u, \na^3 \h) \|_{L^2 (\mbR^3)}$ to control $ \| (\De \om^u, \De \om^h )\|_{L^2 (\mc{C}_{2, \la })} $. These verify \eqref{omthetadecay} and \eqref{omrzdecay}.

\section{Proof of Theorem \ref{main12}}\noindent

In this section, we will consider the special class of axisymmetric D-solution to the steady MHD equations, where the magnetic field has only the swirl component $\textbf{h}(r,z)= h_{\th} (r,z) \bf{e_{\th}}$. In this case, ${\bf \omega^h}(r,z)=\nabla\times \textbf{h(x)}=\omega_r^h(r,z){\bf e_r}+\omega_z^h(r,z){\bf e_z}$, the equations (\ref{vorticity eq}) reduce to
\begin{equation}\lab{special-vorticity}
\left\{
\begin{array}{l}
(u_r\partial_r+u_z\partial_z)\omega_r^u-(\omega_r^u\partial_r+\omega_z^u\partial_z)u_r=(\partial^2_r+\frac{1}{r}\partial_r+\partial^2_z-\frac{1}{r^2})\omega_r^u
- \p_z f_{\th},\\
(u_r\partial_r+u_z\partial_z)\omega_{\theta}^u-\frac{u_r\omega_{\theta}^u}{r}+\frac{1}{r}\partial_z(h_{\theta}^2-u_{\theta}^2)=(\partial^2_r+\frac{1}{r}\partial_r+\partial^2_z-\frac{1}{r^2})\omega_{\theta}^u
+ (\p_z f_r - \p_r f_z),\\
(u_r\partial_r+u_z\partial_z)\omega_z^u-(\omega_r^u\partial_r+\omega_z^u\partial_z)u_z=(\partial^2_r+\frac{1}{r}\partial_r+\partial^2_z)\omega_z^u
+ \f 1r \p_r (r f_{\th}),\\
(u_r\partial_r+u_z\partial_z)\omega_r^h-(\omega_r^h\partial_r+\omega_z^h\partial_z)u_r+\frac{2}{r}\partial_z(u_rh_{\theta})=(\partial^2_r+\frac{1}{r}\partial_r+\partial^2_z-\frac{1}{r^2})\omega_r^h\\
\q\q - \p_z (\p_z g_r - \p_r g_z),\\
(u_r\partial_r+u_z\partial_z)\omega_z^h-(\omega_r^h\partial_r+\omega_z^h\partial_z)u_z-\frac{2}{r}\partial_r(u_rh_{\theta})=(\partial^2_r+\frac{1}{r}\partial_r+\partial^2_z)\omega_z^h
\\ \q\q+ \f 1r \p_r (r (\p_z g_r - \p_r g_z)).
\end{array}
\right.
\end{equation}

Comparing with \eqref{vorticity eq}, \eqref{special-vorticity} has a simpler form, from which we derive better localized energy estimates and it turns out that we can obtain better decay rates.

The proof of theorem \ref{main12} will be divided into the following steps. Firstly, the decay rate of $\om_{\th}^u$ and a weaker decay rates of $(\om_r^u, \om_z^u)$ will be deduced by using the simpler structure of the equations in \eqref{special-vorticity}. Secondly, note that $\na (u_r \e_r + u_z \e_z ) = \mb{C}(\om_{\th}^u \e_{\th}) + \mb{K} \ast (\om_{\th}^u \e_{\th})$, where $\mathbb{C}(\om_{\th}^u \e_{\th})$ is a matrix whose entries are components of $\om_{\th}^u \e_{\th}$ and $\mathbb{K}$ is a \textit{Calderon-Zymund kernel}, then Lemma \ref{lemma CZ} and the estimate \eqref{decay-om-th-u} will be used to improve the decay rates of $(\om_r^u, \om_z^u)$:
  \be\no
  |\na u_r| + |\na u_z| \le C (M) r^{-5/4}(\ln r)^{7/4},
  \ee
  and
\be\no
&& |\om_r^u (r,z)| + |\om_z^u (r,z)| \le C (M ) r^{- 9/8} (\ln r)^{11/8},\\\no
& &
|\om_r^h (r,z)| + |\om_z^h (r,z)| \le C (M) r^{-1} (\ln r)^{11/8},
\ee
for large $r$. Thirdly, we utilize the weighted energy estimates to obtain better decay rates of $\na h_{\th }$:
\be\no
|\om^h_r|+ |\om^h_z| + |\na h_{\th} | \leq C(M) r^{- (\f {35}{32})^-}.
\ee

\textbf{Step 1: }Proof of \eqref{decay-om-th-u} for $\omega_{\theta}^u$ and a weaker decay of $(\omega_r^u, \omega_z^u)$.

Same as in section 3, consider the scaled solution $\tilde{\textbf{u}}(\tilde{\textbf{x}})$, $\tilde{\textbf{h}}(\tilde{\textbf{x}})$, $\tilde{\om}(\tilde{\textbf{x}})$, and $\tilde{\ff} (\tilde{\textbf{x}})$, $\curl_{\tilde{\textbf{x}}} \tilde{\g} (\tilde{\textbf{x}})$, where $\tilde{\textbf{x}}=\lambda {\textbf{x}}$. Drop the "$\sim$" for simplification of notation when computations take place under the scaled sense. Taking the same cut-off function $\phi(y)$ as previous, and testing the vorticity equation (\ref{special-vorticity}) with $\omega_r^u\phi^2$, $\omega_{\theta}^u\phi^2$, $\omega_z^u\phi^2$, $\omega_r^h\phi^2$ and $\omega_z^h\phi^2$ respectively, we deduce that
\be\no
&&\int_{\mathcal{C}_1}\bigg(|\nabla(\omega_{\theta}^u\phi)|^2+\frac{(\omega_{\theta}^u\phi)^2}{r^2}\bigg)dy=\int_{\mathcal{C}_1}\bigg(
|\omega_{\theta}^u|^2|\nabla\phi|^2+\frac{1}{2}|\omega_{\theta}^u|^2 (u_r \p_r + u_z \p_z) (\phi^2)\\\no
&&\quad\quad+\frac{u_r}{r}(\omega_{\theta}^u)^2\phi^2- 2 \frac{\omega_{\theta}^u\phi^2}{r}(-\omega_r^hh_{\theta}+\omega_r^uu_{\theta})
+ f_z \f {\om_{\th}^u \phi}{r} \phi  \\\no
&&\q\q + \phi(f_z \p_r - f_r \p_z) (\om_{\th}^u \phi) + \om_{\th}^u \phi (f_z \p_r - f_r \p_z) \phi
\bigg)dy\\\no
&&\leq C (1+\|(u_r,u_{\theta},u_z,h_{\theta})\|_{L^{\infty}(\mathcal{C}_1)})\|(\omega_r^u,\omega_{\theta}^u,\omega_r^h)\|^2_{L^2(\mathcal{C}_1)}
+ C \|(f_r, f_z) \|_{L^2 (\mc{C}_1)}^2\\\no
&& \q \q
+ \f 18 \| \na (\om_{\th}^u \phi) \|^2_{L^2 (\mathcal{C}_1)},
\ee
\be\no
&&\int_{\mathcal{C}_1}\bigg(|\nabla(\omega_r^u\phi)|^2+\frac{(\omega_r^u\phi)^2}{r^2}\bigg)dy=\int_{\mathcal{C}_1}\bigg(
(\omega_r^u)^2|\nabla\phi|^2+ \f 12 |\om_r^u|^2 (u_r \p_r + u_z \p_z) (\phi^2) \\\no
&&\quad - u_r (\p_r((\om_r^u \phi)^2) + \p_z (\om_r^u \phi \cdot \om_z^u \phi))- \f {(\om_r^u \phi)^2}{r} u_r
+ f_{\th}\phi (\p_z (\om_r^u \phi) + \om_r^u \p_z \phi)
\bigg)dy\\\no
&&=\int_{\mathcal{C}_1}\bigg((\omega_r^u)^2|\nabla\phi|^2+\frac{1}{2}(\omega_r^u)^2(u_r\partial_r+u_z\partial_z)(\phi^2)
+ (|\om_r^u|^2 \phi^2 \p_r + \om_r^u \phi \cdot \om_z^u \phi \p_z) u_r \\\no
&& \q\q + f_{\th}\phi (\p_z (\om_r^u \phi) + \om_r^u \p_z \phi)
\bigg)dy\\\no
&&\leq \begin{cases}
C  (1+\|(u_r,u_z)\|^2_{L^{\infty}(\mathcal{C}_1)})\|(\omega_r^u, \om_z^u)\|^2_{L^2(\mathcal{C}_1)}+ C  \| f_{\th } \|_{L^2 (\mc{C}_1)}^2\\\no
\q\q +\frac{1}{8}\|(\nabla(\omega_r^u\phi), \na(\om_z^u \phi)\|^2_{L^2(\mathcal{C}_1)},\\\no
C  (1+\|(u_r,u_z)\|_{L^{\infty}(\mathcal{C}_1)}+\|\nabla u_r\|_{L^{\infty}(\mathcal{C}_1)})\|(\omega_r^u,\omega_z^u)\|^2_{L^2(\mathcal{C}_1)}
+ C \| f_{\th } \|_{L^2 (\mc{C}_1)}^2\\\no
 \q\q +\frac{1}{8}\|\nabla(\omega_r^u\phi)\|^2_{L^2(\mathcal{C}_1)},
\end{cases}\ee
\be\no
&&\int_{\mathcal{C}_1}|\nabla(\omega_z^u\phi)|^2dy=\int_{\mathcal{C}_1}\bigg(
(\omega_z^u)^2|\nabla\phi|^2+\frac{1}{2}(\omega_z^u)^2(u_r\partial_r+u_z\partial_z)(\phi^2)\\\no
&&\quad\quad- u_z \p_r(\om_r^u \phi \cdot \om_z^u \phi) - u_z \p_z(\om_z^u \phi)^2- \f {u_z}{r} \om_r^u \om_z^u \phi^2
- f_{\th} \phi (\p_r (\om_z^u \phi) + \om_z^u \p_r \phi)
\bigg)dy\\\no
&&=\int_{\mathcal{C}_1}\bigg((\omega_z^u)^2|\nabla\phi|^2+\frac{1}{2}(\omega_z^u)^2(u_r\partial_r+u_z\partial_z)(\phi^2)
+ (\omega_z^u)^2\phi^2\partial_z u_z\\\no
&&\q\q+ \omega_r^u\omega_z^u\phi^2\partial_r u_z  - f_{\th} \phi (\p_r (\om_z^u \phi) + \om_z^u \p_r \phi)
\bigg)dy\\\no
&&\leq \begin{cases}
C  (1+\|(u_r,u_z)\|^2_{L^{\infty}(\mathcal{C}_1)})\|(\om_r^u,\omega_z^u)\|^2_{L^2(\mathcal{C}_1)}
+ C \|f_{\th } \|_{L^2 (\mc{C}_1)}^2 \\\no
\q\q +\frac{1}{8} \|(\nabla(\omega_r^u\phi),\nabla(\omega_z^u\phi)) \|^2_{L^2(\mathcal{C}_1)}\\\no
C  (1+\|(u_r,u_z)\|_{L^{\infty}(\mathcal{C}_1)}+\|\nabla u_z\|_{L^{\infty}(\mathcal{C}_1)})\|(\omega_r^u,\omega_z^u)\|^2_{L^2(\mathcal{C}_1)} + C \| f_{\th } \|_{L^2 (\mc{C}_1)}^2\\\no
\q\q +\frac{1}{8} \|\nabla(\omega_z^u\phi) \|^2_{L^2(\mathcal{C}_1)},
\end{cases}\ee

\be\no
&&\int_{\mathcal{C}_1}\bigg(|\nabla(\omega_r^h\phi)|^2+\frac{(\omega_r^h\phi)^2}{r^2}\bigg)dy=
\int_{\mathcal{C}_1}\bigg(|\omega_r^h|^2|\nabla\phi|^2+\frac{1}{2}(\omega_r^h)^2(u_r\partial_r+u_z\partial_z)(\phi^2)
\\\no
&&\quad\quad- u_r \p_r(\om_r^h \phi)^2 - u_r \p_z (\om_r^h \phi \cdot \om_z^h \phi) - \f {u_r}{r} (\om_r^h \phi)^2 + 2 \f {u_r h_{\th}}{r} \p_z (\om_r^h \phi^2)\\\no
&& \q\q + \p_z (\om_r^h \phi) (\p_z g_r - \p_r g_z)
\bigg)dy\\\no
&&= \int_{\mathcal{C}_1}\bigg(|\omega_r^h|^2|\nabla\phi|^2
+\frac{1}{2}|\omega_r^h|^2(u_r\partial_r+u_z\partial_z)(\phi^2)
+ ((\om_r^h)^2 \phi^2 \p_r + \om_r^h \om_z^h \phi^2 \p_z ) u_r\\\no
&& \q\q - 2 (\p_r u_r + \p_z u_z) h_{\th} \p_z(\om_r^h \phi^2)
 + \p_z (\om_r^h \phi) (\p_z g_r - \p_r g_z)
\bigg)dy\\\no
&&\leq \begin{cases}
C  (1+\|(u_r,u_z)\|^2_{L^{\infty}(\mathcal{C}_1)})\|(\omega_r^h,\omega_z^h)\|^2_{L^2(\mathcal{C}_1)}
+\frac{1}{8}\|(\nabla(\omega_r^h\phi),\nabla(\omega_z^h\phi))\|^2_{L^2(\mathcal{C}_1)}\\
\quad\quad+ C (M) \|u_r\|^2_{L^{\infty}(\mathcal{C}_1)}\int_{\mbR^3} \f {h_{\th}^2}{r^2} dy + C \| \bar {\g} \|_{L^2 (\mc{C}_1)}^2 \\
C  (1+\|(u_r,u_z)\|_{L^{\infty}(\mathcal{C}_1)}+\|\nabla u_r\|_{L^{\infty}(\mathcal{C}_1)})\|(\omega_r^h,\omega_z^h)\|^2_{L^2(\mathcal{C}_1)}
+\frac{1}{8}\|\nabla(\omega_r^h\phi)\|^2_{L^2(\mathcal{C}_1)}
\\
\quad\quad+C  \|h_{\theta}\|^2_{L^{\infty}(\mathcal{C}_1)}\|(\nabla u_r, \na u_z)\|^2_{L^2(\mathcal{C}_1)} + C \| \bar {\g} \|_{L^2 (\mc{C}_1)}^2 ,
\end{cases}\ee
\be\no
&&\int_{\mathcal{C}_1}|\nabla(\omega_z^h\phi)|^2dy=\int_{\mathcal{C}_1}\bigg(|\omega_z^h|^2|\nabla\phi|^2
+\frac{1}{2}|\omega_z^h|^2(u_r\partial_r+u_z\partial_z)(\phi^2)- u_z \p_r (\om_r^h \phi \cdot \om_z^h \phi) \\\no
&&\quad\quad - \f {u_z}{r} \om_r^h \om_z^h \phi^2 - u_z \p_z (\om_z^h \phi)^2 - 2 \f {u_r h_{\th}}{r} \p_r (\om_z^h \phi^2)
- \p_r (\om_z^h \phi^2) (\p_z g_r - \p_r g_z)
\bigg)dy\\\no
&&=\int_{\mathcal{C}_1}\bigg(|\omega_z^h|^2|\nabla\phi|^2+\frac{1}{2}|\omega_z^h|^2(u_r\partial_r+u_z\partial_z)(\phi^2)
+ (\om_r^h \om_z^h \phi^2 \p_r + (\om_z^h \phi)^2 \p_z) u_z\\\no
&& \quad\quad+ 2 \f {\om_z^h \phi^2} r \p_r (u_r h_{\th})- \p_r (\om_z^h \phi^2) (\p_z g_r - \p_r g_z)
\bigg) dy\\\no
&&\leq \begin{cases}
C  (1+\|(u_r,u_z)\|_{L^{\infty}(\mathcal{C}_1)}^2)\|(\omega_r^h,\omega_z^h)\|^2_{L^2(\mathcal{C}_1)}
+\frac{1}{8}\|(\nabla(\omega_r^h\phi),\nabla(\omega_z^h\phi))\|^2_{L^2(\mathcal{C}_1)}\\
\quad\quad+C
\|u_r\|^2_{L^{\infty}(\mathcal{C}_1)}\int_{\mbR^3} \f {h_{\th}^2}{r^2} dy + C \| \bar {\g} \|_{L^2 (\mc{C}_1)}^2 ,\\
C  (1+\|(u_r,u_z)\|_{L^{\infty}(\mathcal{C}_1)}+\|\nabla u_z\|_{L^{\infty}(\mathcal{C}_1)})\|(\omega_r^h,\omega_z^h)\|^2_{L^2(\mathcal{C}_1)}
 +\frac{1}{8}\|\nabla(\omega_z^h\phi)\|^2_{L^2(\mathcal{C}_1)}
\\
\quad\quad+C
\|(\nabla u_r, \na h_{\th})\|^2_{L^2(\mathcal{C}_1)}\|( u_r, h_{\theta})\|^2_{L^{\infty}(\mathcal{C}_1)} + C \| \bar {\g} \|_{L^2 (\mc{C}_1)}^2 .
\end{cases}\ee

Thus we deduce
\be\lab{om-th}
&& \|\nabla\omega_{\theta}^u\|^2_{L^2(\mathcal{C}_2)} \leq  C(1+\|(u_r,u_{\theta},u_z,h_{\theta})\|_{L^{\infty}(\mathcal{C}_1)})\|(\omega_r^u,\omega_{\theta}^u,\omega_r^h)\|^2_{L^2(\mathcal{C}_1)}\\\no
&& \q \q \q\q\q\q\q
 + C \|(f_r, f_z) \|_{L^2 (\mc{C}_1)}^2,
\\\lab{om-rzu}
&& \|(\nabla\omega_r^u,\nabla\omega_z^u)\|^2_{L^2(\mathcal{C}_2)} \leq C(1+\|(u_r,u_z)\|^2_{L^{\infty}(\mathcal{C}_1)})\|(\omega_r^u,\omega_z^u)\|^2_{L^2(\mathcal{C}_1)}
+ C \| f_{\th } \|_{L^2 (\mc{C}_1)}^2,
\\\lab{om-rzh}
&& \|(\nabla\omega_r^h,\nabla\omega_z^h)\|^2_{L^2(\mathcal{C}_2)} \leq  C(1+\|(u_r,u_z)\|^2_{L^{\infty}(\mathcal{C}_1)})\|(\omega_r^h,\omega_z^h)\|^2_{L^2(\mathcal{C}_1)} \\\no
&& \q \q \q \q \q \q \q\q\q\q+ C \| u_r \|^2_{L^{\oo} (\mathcal{C}_1)} +C  \| \bar {\g} \|_{L^2 (\mc{C}_1)}^2,
\ee
and
\be\no
&& \|\nabla(\omega_r^u,\omega_z^u)\|^2_{L^2(\mathcal{C}_2)} \leq  C(1+\|(u_r,u_z)\|_{L^{\infty}(\mathcal{C}_1)}+\|\nabla (u_r,u_z)\|_{L^{\infty}(\mathcal{C}_1)})\|(\omega_r^u,\omega_z^u)\|^2_{L^2(\mathcal{C}_1)}
\\\label{na-om-rzu}
&& \q\q\q\q\q\q\q \q \q + C \| f_{\th } \|_{L^2 (\mc{C}_1)}^2,\\\no
&& \|\nabla(\omega_r^h,\omega_z^h)\|^2_{L^2(\mathcal{C}_2)} \leq  C(1+\|(u_r,u_z)\|_{L^{\infty}(\mathcal{C}_1)}+\|(\nabla u_r,\nabla u_z)\|_{L^{\infty}(\mathcal{C}_1)})\|(\omega_r^h,\omega_z^h)\|^2_{L^2(\mathcal{C}_1)}\\\lab{na-om-rzh}
&& \q\q\q\q\q\q\q \q \q
+C\|(\nabla u_r,\na h_{\th})\|^2_{L^2(\mathcal{C}_1)}\|(u_r, h_{\theta})\|^2_{L^{\infty}(\mathcal{C}_1)}
 + C \| \bar {\g} \|_{L^2 (\mc{C}_1)}^2.
\ee

As previous, utilizing the \textit{Brezis-Gallouet's inequality} and (\ref{om-th})-(\ref{na-om-rzh}), we are led to
\be\no
&&\|\omega_{\theta}^u\|_{L^{\infty}(\mathcal{\bar{C}}_2)}\leq C\bigg(1+(1+\|(u_r,u_{\theta},u_z,h_{\theta})\|^{1/2}_{L^{\infty}(\mathcal{C}_1)})\|(\omega_r^u,\omega_{\theta}^u,\omega_r^h)\|_{L^2(\mathcal{C}_1)}
+ \| (f_r, f_z) \|_{L^2 (\mc{C}_1)} \bigg)\\\no
&&\quad\quad\quad \quad\quad\times \ln^{1/2}\bigg(e+\|\De\omega_{\theta}^u\|_{L^2(\mathcal{\bar{C}}_2)}\bigg),
\ee
\be\no
&&\|(\omega_r^u,\omega_z^u)\|_{L^{\infty}(\mathcal{\bar{C}}_2)}\leq C\bigg(1+(1+\|(u_r,u_z)\|_{L^{\infty}(\mathcal{C}_1)})\|(\omega_r^u,\omega_z^u)\|_{L^2(\mathcal{C}_1)}
+ \| f_{\th } \|_{L^2 (\mc{C}_1)}  \bigg)\\\no
&&\quad\quad\quad \quad\quad\times\ln^{1/2}\bigg(e+\|(\De\omega_r^u,\De\omega_z^u)\|_{L^2(\mathcal{\bar{C}}_2)}\bigg),\\\no
&&\|(\omega_r^h,\omega_z^h)\|_{L^{\infty}(\mathcal{\bar{C}}_2)}\leq C\bigg(1+(1+\|(u_r,u_z)\|_{L^{\infty}(\mathcal{C}_1)})\|(\omega_r^h,\omega_z^h)\|_{L^2(\mathcal{C}_1)} + \| u_r \|_{L^{\infty}(\mathcal{C}_1)}\\\no
&& \quad\quad\quad \quad\quad
+ \| \bar {\g} \|_{L^2 (\mc{C}_1)}
\bigg) \ln^{1/2}\bigg(e+\|(\De\omega_r^u,\De\omega_z^u)\|_{L^2(\mathcal{\bar{C}}_2)}\bigg).
\ee
and
\be\no
&&\|(\omega_r^u,\omega_z^u)\|_{L^{\infty}(\mathcal{\bar{C}}_2)}\leq C\bigg(1+(1+\|(u_r,u_z)\|^{1/2}_{L^{\infty}(\mathcal{C}_1)}+\|(\nabla u_r,\nabla u_z)\|^{1/2}_{L^{\infty}(\mathcal{C}_1)})\|(\omega_r^u,\omega_z^u)\|_{L^2(\mathcal{C}_1)}\\\no
&& \quad \quad + \| f_{\th } \|_{L^2 (\mc{C}_1)} \bigg) \ln^{1/2}\bigg(e+\|(\De\omega_r^u,\De\omega_z^u)\|_{L^2(\mathcal{\bar{C}}_2)}\bigg),\\\no
&&\|(\omega_r^h,\omega_z^h)\|_{L^{\infty}(\mathcal{\bar{C}}_2)}\leq C\bigg(1+(1+\|(u_r,u_z)\|^{1/2}_{L^{\infty}(\mathcal{C}_1)}+\|(\nabla u_r,\nabla u_z)\|^{1/2}_{L^{\infty}(\mathcal{C}_1)})\|(\omega_r^h,\omega_z^h)\|_{L^2(\mathcal{C}_1)}\\\no
&& \quad\quad +\|(\nabla u_r,\na h_{\th})\|_{L^2(\mathcal{C}_1)}
\|(u_r,h_{\theta})\|_{L^{\infty}(\mathcal{C}_1)}
+ \|\bar {\g} \|_{L^2 (\mc{C}_1)} \bigg)
\ln^{1/2}\bigg(e+\|(\De\omega_r^u,\De\omega_z^u)\|_{L^2(\mathcal{\bar{C}}_2)}\bigg).
\ee
Then scaling back to find
\be\no
&&\lambda^2\|\omega_{\theta}^u\|_{L^{\infty}(\mathcal{C}_{2,\lambda})}\leq C\bigg(1+(1+\lambda^{1/2} \|(u_r,u_{\theta},u_z,h_{\theta})\|^{1/2}_{L^{\infty}(\mathcal{C}_{1,\lambda})})\lambda^{1/2} \|(\omega_r^u,\omega_{\theta}^u,\omega_r^h)\|_{L^2(\mathcal{C}_{1,\lambda})}
\\\no
&& \quad \quad + \la^{3/2} \| (f_r, f_z) \|_{L^2 (\mc{C}_{1, \la })} \bigg) \ln^{1/2}\bigg(e+\lambda^{5/2}\|\De\omega_{\theta}^u\|_{L^2(\mathcal{C}_{2,\lambda})}\bigg),
\\\no
&&\lambda^2\|(\omega_r^u,\omega_z^u)\|_{L^{\infty}(\mathcal{C}_{2,\lambda})}\leq C\lambda^{1/2}\bigg(1+(1+\lambda\|(u_r,u_z)\|_{L^{\infty}(\mathcal{C}_{1,\lambda})})\|(\omega_r^u,\omega_z^u)\|_{L^2(\mathcal{C}_{1,\lambda})}
\\\no
&& \quad \quad
+ \la^{3/2} \| f_{\th } \|_{L^2 (\mc{C}_{1, \la })} \bigg) \ln^{1/2}\bigg(e+\lambda^{5/2}\|(\De\omega_r^u,\De\omega_z^u)\|_{L^2(\mathcal{C}_{2,\lambda})}\bigg),\\\no
&&\lambda^2\|(\omega_r^h,\omega_z^h)\|_{L^{\infty}(\mathcal{C}_{2, \la})}\leq C\bigg(1+(1+\lambda\|(u_r,u_z)\|_{L^{\infty}(\mathcal{C}_{1,\lambda})})\lambda^{1/2}\|(\omega_r^h,\omega_z^h)\|_{L^2(\mathcal{C}_{1,\lambda})}
\\\no
&&\quad \quad
+\lambda\|u_r\|_{L^{\infty}(\mathcal{C}_{1,\lambda})}
+ \la^{ 3/2} \|\bar { \g} \|_{L^2 (\mc{C}_{1, \la })}  \bigg)
\ln^{1/2}\bigg(e+\lambda^{5/2}\|(\De\omega_r^h,\De\omega_z^h)\|_{L^2(\mathcal{C}_{2,\lambda})}\bigg),
\ee
and
\be\lab{na-om-rzu-lam}
&&\lambda^2\|(\omega_r^u,\omega_z^u)\|_{L^{\infty}(\mathcal{C}_{2,\lambda})}\le  C\lambda^{1/2}\bigg(1+(1+\lambda^{1/2}\|(u_r,u_z)\|^{1/2}_{L^{\infty}(\mathcal{C}_{1,\lambda})}
+ \la^{3/2} \| f_{\th } \|_{L^2 (\mc{C}_{1, \la })} \\\no
&&\q +\lambda\|(\nabla u_r,\nabla u_z)\|^{1/2}_{L^{\infty}(\mathcal{C}_{1,\lambda})})\|(\omega_r^u,\omega_z^u)\|_{L^2(\mathcal{C}_{1,\lambda})}\bigg) \ln^{1/2}\bigg(e+\lambda^{5/2}\|(\De\omega_r^u,\De\omega_z^u)\|_{L^2(\mathcal{C}_{2,\lambda})}\bigg),
\ee
\be\no
&&\lambda^2\|(\omega_r^h,\omega_z^h)\|_{L^{\infty}(\mathcal{C}_{2, \la})}
\leq  C \bigg(
1+(1+\lambda^{1/2}\|(u_r,u_z)\|^{1/2}_{L^{\infty}(\mathcal{C}_{1,\lambda})}
+\lambda \|(\nabla u_r,\nabla u_z)\|^{1/2}_{L^{\infty}(\mathcal{C}_{1,\lambda})})\\\no
&& \q
 \times \lambda^{1/2} \|(\omega_r^h,\omega_z^h)\|_{L^2(\mathcal{C}_{1,\lambda})} + \lambda^{\frac32}\|(\nabla u_r,\na h_{\th})\|_{L^2(\mathcal{C}_{1,\lambda})} \|(u_r, h_{\theta})\|_{L^{\infty}(\mathcal{C}_{1,\lambda})}
+ \la^{ 3/2} \| \bar {\g} \|_{L^2 (\mc{C}_{1, \la })}  \bigg)\\\lab{na-om-rzh-lam}
&&\quad \times
\ln^{1/2}\bigg(e+\lambda^{5/2}\|(\De\omega_r^h,\De\omega_z^h)\|_{L^2(\mathcal{\bar{C}}_2)}\bigg).
\ee
 By the a priori bound in \eqref{uhdecay}, we have
\be\no
\|\omega_{\theta}^u\|_{L^{\infty}(\mathcal{C}_{2,\lambda})} & \le &  C (M) \f {(\ln \la)^{3/4}}{\la^{5/4}}  ,\\
\|(\omega_r^u,\omega_z^u)\|_{L^{\infty}(\mathcal{C}_{2,\lambda})}& \le &  C (M) \frac{\ln \lambda}{\lambda},\\\no
\|(\omega_r^h,\omega_z^h)\|_{L^{\infty}(\mathcal{C}_{2,\lambda})}& \le &  C (M)  \frac{\ln \lambda}{\lambda}.
\ee
This verifies (\ref{decay-om-th-u}).

We further derive the decay of $\nabla \omega_{\theta}^u$ under the additional assumption that $\ff \in H^2(\mathbb{R}^3)$ and $ \g \in H^3(\mathbb{R}^3)$. Recall that
 \begin{equation}\lab{4-1}
 (u_r\partial_r+u_z\partial_z)\omega_{\theta}^u-\frac{u_r\omega_{\theta}^u}{r}+\frac{1}{r}\partial_z(h_{\theta}^2-u_{\theta}^2)=(\partial^2_r+\frac{1}{r}\partial_r+\partial^2_z-\frac{1}{r^2})\omega_{\theta}^u
+ (\p_z f_r - \p_r f_z).
 \end{equation}
Choose the domain
 \be\no
\mc{C}_3 = \{ (r, \th, z): \f 78 < r < \f 98, 0 \le \th \le 2 \pi, |z| \le \f 14 \}.
 \ee
Take a cut-off function $\phi (y) \in C_0^{\oo} (\mc{C}_2)$, satisfying $0 \le \phi \le 1 $, $\phi (y) = 1 $ for $y \in \mc{C}_3$ and $ | \na \phi | \le C $, for some positive constant $C$. Taking $\p_r$ to the equation in \eqref{4-1} and multiplying the resulting equation by $\p_r \om_{\th}^u \phi^2$, and integrating over $\mathcal{C}_2$, we deduce that
  \be\no
&& \int_{\mc{C}_2} |\na (\p_r \om_{\th}^u \phi)|^2 dy
+ \int_{\mc{C}_2} 2 \f{|\p_r \om_{\th}^u|^2 \phi^2}{r^2} dy
\\\no
&&
=  \int_{\mc{C}_2} |\p_r \om_{\th}^u|^2 |\na \phi|^2 dy
+ \int_{\mc{C}_2} \f 2{r^3} \om_{\th}^u \p_r \om_{\th}^u \phi^2 dy
+ \int_{\mc{C}_2} u_z \f 1r \p_z \om_{\th}^u \p_r \om_{\th}^u \phi^2 dy\\\no
&& +
\int_{\mc{C}_2} 2 u_r  \p_r (\p_r \om_{\th}^u \phi) \p_r \om_{\th}^u \phi dy
- \int_{\mc{C}_2} \f 1{r^2} u_r \om_{\th}^u \p_r \om_{\th}^u \phi^2 dy
+ \int_{\mc{C}_2} \f 1r u_r |\p_r \om_{\th}^u|^2 \phi^2 dy
\\\no
&&+ \int_{\mc{C}_2} u_z [\p_z (\p_r \om_{\th}^u \phi) \p_r \om_{\th}^u \phi - | \p_r \om_{\th}^u|^2 \phi \p_z \phi + \p_z \om_{\th}^u \p_r (\p_r \om_{\th}^u \phi) \phi + \p_z \om_{\th}^u \p_r \om_{\th}^u \phi \p_r \phi ] dy\\\no
&&
- \int_{\mc{C}_2} (u_r \p_r + u_z \p_z) (\p_r \om_{\th}^u \phi) \p_r \om_{\th}^u \phi - (u_r \p_r + u_z \p_z) \phi |\p_r \om_{\th}^u|^2 \phi dy\\\no
&&
- \int_{\mc{C}_2} u_r \bigg[ \f 1r \om_{\th}^u \p_r (\p_r \om_{\th}^u \phi ) \phi + \f 1r \om_{\th}^u \p_r \om_{\th}^u \phi \p_r \phi \bigg]dy
+ \int_{\mc{C}_2} \f 2{r^2} (u_{\th} \om_r^u  - h_{\th} \om_r^h ) \p_r \om_{\th}^u \phi^2  dy\\\no
&& +
2 \int_{\mc{C}_2} h_{\th} [\p_r (\p_r \om_{\th}^u \phi ) \f {\om_r^h \phi}{r} + \f 1r \p_r \om_{\th}^u \phi \om_r^h \p_r \phi ]
-
u_{\th } [\p_r (\p_r \om_{\th}^u \phi ) \f {\om_r^u \phi}r + \f 1r \p_r \om_{\th}^u \phi \om_r^u \p_r \phi] dy\\\no
&& - \int_{\mc{C}_2} (\p_z f_r - \p_r f_z) (\p_r (\p_r \om_{\th}^u \phi) \phi  + \p_r \om_{\th}^u \phi \p_r \phi + \f 1r \p_r \om_{\th}^u \phi^2) dy
\\\no
& \le & C  (1 + \| (\u, \textbf{h} ) \|_{L^{\oo} (\mc{C}_2)} + \| (\u, \textbf{h} ) \|_{L^{\oo} (\mc{C}_2)}^2 )
\| \na \om_{\th}^u \|_{L^2 (\mc{C}_2)}^2 + C \| \bar{\ff} \|_{L^2 (\mc{C}_2)}^2
\\\no
&&
 +C (1 + \| (\u, \textbf{h} ) \|_{L^{\oo} (\mc{C}_2)}^2 ) \| (\om_r^u,\om_{\th}^u, \om_r^h) \|^2_{L^2 (\mc{C}_2)} + \f 18 \| \na (\p_r \om_{\th}^u \phi)\|_{L^2 (\mc{C}_2)}^2\\\no
& \leq &
 C  (1+\|(\u, \textbf{h} )\|_{L^{\infty}(\mathcal{C}_2)} + \|(\u, \textbf{h} )\|_{L^{\infty}(\mathcal{C}_2)}^2 + \|(\u, \textbf{h} )\|_{L^{\infty}(\mathcal{C}_2)}^3) \| (\om_r^u,\om_{\th}^u, \om_r^h) \|^2_{L^2 (\mc{C}_2)}\\\no
  && + (1 + \| (\u, \textbf{h} ) \|_{L^{\oo} (\mc{C}_2)} + \| (\u, \textbf{h} ) \|_{L^{\oo} (\mc{C}_2)}^2 ) \| (f_r, f_z) \|_{L^2 (\mc{C}_2)}^2 + C \| \bar{\ff} \|_{L^2 (\mc{C}_2)}^2\\\no
   && + \f 18 \| \na (\p_r \om_{\th}^u \phi)\|_{L^2 (\mc{C}_2)}^2,
\ee
where the second inequality follows from \eqref{om-th}.

Set
\be\no
\bar{\mc{C}}_3 \= \{ (r,z): \f 78 < r< \f 98, |z| \le \f 14 \}.
\ee
Utilizing the \textit{Brezis-Gallouet} inequality in Lemma \ref{BG ineq} and the above localized energy estimates, we conclude
\be\no
&& \| \p_r \om_{\th}^u \|_{L^{\oo} (\bar{\mc{C}}_3)}  \\\no
& \leq & C
(1 + \| \na (\p_r \om_{\th}^u) \|_{L^2 (\mc{C}_3)}) \ln^{1/2} (e + \| \De (\p_r \om_{\th}^u )\|_{L^2 (\mc{C}_3)})\\\no
&\leq &
C  \bigg(1 +
(1+\|(\u, \textbf{h} )\|_{L^{\infty}(\mathcal{C}_2)}^{\f 12} + \|(\u, \textbf{h} )\|_{L^{\infty}(\mathcal{C}_2)} + \|(\u, \textbf{h} )\|_{L^{\infty}(\mathcal{C}_2)}^{\f 32}) \| (\om_r^u,\om_{\th}^u, \om_r^h) \|_{L^2 (\mc{C}_2)} \\\no
&& \q   + (1 + \| (\u, \textbf{h} ) \|^{\f 12}_{L^{\oo} (\mc{C}_2)} + \| (\u, \textbf{h} ) \|_{L^{\oo} (\mc{C}_2)} )\| (f_r, f_z ) \|_{L^2 (\mc{C}_2)} + \| \bar{\ff} \|_{L^2 (\mc{C}_2)}\bigg)\\\no
&& \q
\times
  \ln^{1/2} (e + \| \De (\p_r \om_{\th}^u )\|_{L^2 (\mc{C}_3)}).
\ee
Then scaling back to the domains $\mc{C}_{2,\la}$ and
\be\no
\mc{C}_{3,\la} = \{ (r, \th, z): \f {7 \la}{8} < r < \f {9 \la}{8}, 0 \le \th \le 2 \pi, |z| \le \f {\la }{4} \},
\ee
we derive that
\be\no
& &\la^3 \| \p_r \om_{\th}^u \|_{L^{\oo} (\mc{C}_{3, \la})} \\\no
& \leq &
C  \bigg( 1 + ( 1 + \la^{\f 12} \|(\u, \textbf{h} )\|_{L^{\infty}(\mathcal{C}_{2, \la })}^{\f 12}+ \la \|(\u, \textbf{h} )\|_{L^{\infty}(\mathcal{C}_{2, \la})} + \la^{\f 32}\|(\u, \textbf{h} )\|_{L^{\infty}(\mathcal{C}_{2, \la})}^{\f 32})
 \\\no
&& \q \times \la^{\f 12}\| (\om_r^u,\om_{\th}^u, \om_r^h) \|_{L^2 (\mc{C}_{2, \la})}
+ (1 + \la^{\f 12} \|(\u, \textbf{h} )\|_{L^{\infty}(\mathcal{C}_{2, \la })}^{\f 12}+ \la \|(\u, \textbf{h} )\|_{L^{\infty}(\mathcal{C}_{2, \la})}) \\\no
&& \q \times
\la^{ \f 32} \| (f_r, f_z ) \|_{L^2 (\mc{C}_{2, \la })} + \la^{\f 52} \| \bar{\ff} \|_{L^2 (\mc{C}_{2, \la })}\bigg)
\ln^{\f 12}
(
e + \la^{\f 72} \| \De (\p_r \om_{\th}^u) \|_{L^2 (\mc{C}_{2, \la})}
)\\\no
& \leq & C (M) \la^{\f 32} (\f {\ln \la }{\la })^{\f 34} \la^{\f 12} (\ln \la)^{\f 12}\\\no
& \le & C (M) \la^{- \f 54} (\ln \la)^{\f 54},
\ee
where we use the boundedness of $\|\nabla^4 {\bf u}\|_{L^2(\mathbb{R}^3)}$ to control $\| \De (\p_r \om_{\th}^u) \|_{L^2 (\mc{C}_{2, \la})}$, which follows from the $L^p$ theory to steady Stokes system and a bootstrap argument under the assumption that ${\bf f}\in H^2(\mathbb{R}^3)$ and ${\bf g}\in H^3(\mathbb{R}^3)$.

Thus
\be\no
\| \p_r \om_{\th}^u \|_{L^{\oo} (\mc{C}_{3, \la})}  \le C (M) \f{(\ln \la)^{5/4}}{\la^{ 7/4} }.
\ee
Similarly, there holds
\be\no
\| \p_z \om_{\th}^u \|_{L^{\oo} (\mc{C}_{3, \la})} \le C (M) \f{(\ln \la)^{5/4}}{\la^{ 7/4} }.
\ee
Then
\be\no
\| \na \om_{\th}^u \|_{L^{\oo} (\mc{C}_{3, \la})} \le C (M)\f{(\ln \la)^{5/4}}{\la^{ 7/4} }.
\ee


\textbf{Step 2:} we employ the \textit{Biot-Savart law} to get better decay estimates of $|\omega_r^u|+|\omega_z^u|$ and $|\omega_r^h|+|\omega_z^h|$.


Since $- \De (u_r \e_r + u_z \e_z ) = \curl (\om_{\th }^u \e_{\th })$, then according to the formula (2.105) in \cite[Page 77]{mb},
\be\lab{cz-kernel}
\na (u_r \e_r + u_z \e_z) = \mb{K} \ast (\om_{\th}^u \e_{\th }) + \mathbb{C}(\om_{\th}^u \e_{\th}),
\ee
where $\mathbb{C}(\om_{\th}^u \e_{\th})$ is a matrix whose entries are components of $\om_{\th}^u \e_{\th}$ and $\mathbb{K}$ is a \textit{Calderon-Zymund kernel}.

Applying Lemma \ref{lemma CZ} to (\ref{cz-kernel}) and note the decay of $\om_{\th}^u$, we can get
\begin{equation}\lab{na-u-rz}
|\nabla u_r|+|\nabla u_z|\le C\frac{(\ln \lambda)^{7/4}}{\lambda^{5/4}}.
\end{equation}
Now go back to (\ref{na-om-rzu-lam}) and (\ref{na-om-rzh-lam}), we get better decay on $|\omega_r^u|+ |\omega_z^u|$ and $|\omega_r^h| + |\omega_z^h|$,
\be\no
|\omega_r^u|+|\omega_z^u|&\le& C\lambda^{-3/2}\bigg(\lambda^{1/2}(\lambda^{-1/2}(\ln \lambda)^{1/2})^{1/2}+\lambda(\lambda^{-5/4}(\ln \lambda)^{7/4})^{1/2}\bigg)(\ln \lambda)^{1/2}\\
&\le& C \lambda^{-9/8}(\ln\lambda)^{11/8},\q  \text{ for large $\lambda$},\\\no
|\omega_r^h|+|\omega_z^h|
&\le&
 C\lambda^{-2} \la (\ln \la )^{7/8} (\ln \la)^{1/2}\\
&\le& C \la^{-1} (\ln \la)^{11/8},\q \text{ for large $\lambda$}.
\ee

{\bf Step 3:} Improved decay rate of $\na h_{\th}$ by weighted energy estimates.
By slightly modifying the proof in Lemma 3.2 in \cite{weng15} (taking $\delta=(\frac{1}{2})^-$ in that lemma), we obtain the following estimates, which will be used in the proof of Theorem \ref{main13}.
\bl\lab{weighted estimate1}
{\it Let $({\bf u}, {\bf h}, p)$ be an axially symmetric smooth D-solutions to inhomogeneous stationary MHD equations with \emph{\textbf{f}} and \emph{\textbf{g}} satisfying (\ref{fg}) and (\ref{fg1}), where ${\bf h}(r,z) = h_{\theta} (r,z){\bf e}_{\theta}$. Then the following estimates hold
\begin{eqnarray}\no
&&\int_{\mbR^3}|\om^u_{\th}|^2 dx+\int_{\mbR^3}r^{(\frac{3}{2})^-}|\na\om^u_{\th}|^2 dx + \int_{\mbR^3}r^{(\frac{5}{2})^-}|\p_z\na\om^u_{\th}|^2 dx\leq C(M).
\end{eqnarray}
}
\el

\bl\lab{weighted estimate3}
{\it Let $({\bf u},{\bf h},p)$ be an axially symmetric smooth D-solutions to inhomogeneous stationary MHD equations with \emph{\bf{f}} satisfying (\ref{fg}), (\ref{fg1}) and \emph{\bf{g}} satisfying (\ref{fg2}), where ${\bf h}(r,z) = h_{\th} (r,z){\bf e}_{\th}$, suppose that
\be\lab{uh}
|u_r(r,z)| + |u_{\th}(r,z)| + |u_z(r,z)| + |h_{\th}(r,z)| &\le& C (1+r)^{- \de},\\\lab{na-u}
|\na u_r(r,z)| +|\na u_z(r,z)| &\le& C (1+r)^{-1-\ga},
\ee
holds for some $\de, \ga \in [0,1]$. Then the following estimates holds
\be\lab{naomu}
&&\int_{\mathbb R^3} r^{1+\de \land \ga} (|\na \om_r^u|^2 + |\na \om_z^u|^2) dx \le C(M),\\\lab{naomh}
&&\int_{\mathbb R^3} r^{1+\de \wedge \ga} (|\na \om_r^h|^2 + |\na \om_z^h|^2) dx \le C(M),\\\lab{pnaomu}
&&\int_{\mathbb R^3} r^{1 + \de \wedge \ga + 2\de} (|\partial_z \na \om_r^u|^2 + |\partial_z \na \om_z^u|^2) dx \le C(M),\\\lab{pnaomh}
&&\int_{\mathbb R^3} r^{1 + \de \wedge \ga + 2\de} (|\p_z \na \om_r^h|^2 + |\p_z \na \om_z^h|^2) dx \le C(M),
\ee
where $\de \wedge \ga = \mathrm{min} \{\de, \ga\}$.

In particular, by Lemma \ref{a decay lemma}, we obtain the following decay rate:
\be\no
|\om_r^u(r,z)|+|\om^u_z(r,z)| &\le& C(M) r^{-\frac{7}{8}-\frac {1}{8} (3(\de \wedge \ga) + 2 \de)},\\\no
|\om_r^h(r,z)|+|\om^h_z(r,z)| &\le& C(M) r^{-\frac{7}{8}-\frac {1}{8} (3(\de \wedge \ga) + 2 \de)}.
\ee
}\el
\begin{proof}
According to Lemma 3.8 in \cite{weng15}, we have (\ref{naomu}) and (\ref{pnaomu}). Similarly, we can prove (\ref{naomh}). It remains to prove (\ref{pnaomh}).

We start to derive the weighted estimates of $\na \p_{rz}^2 h_{\th}$.
\be\lab{mag111}
&& \p_{rz}^2 \b[(u_r\p_r+ u_z\p_z) h_{\theta}-\f{u_r}{r} h_{\theta}\b]= \b(\p_r^2+\f1r\p_r+\p_z^2-\f1{r^2}\b)\p_{rz}^2 h_{\theta} \\\no
&&\q\q -\f1{r^2}\p_{rz}^2 h_{\theta}+ \f{2}{r^3}\p_z h_{\theta} + \p_{rz}^2 (\p_z g_r - \p_r g_z)
\ee
take $\eta^2 r^{a_2} \p_{rz}^2 h_{\th}$ as a test function to (\ref{mag111}) and integrating over $\mbR^3$, we obtain
\be\no
0&=& \int_{\mbR^3} \eta^2 r^{a_2} |\na \p_{rz}^2 h_{\th}|^2 dx
+\f12\int_{\mbR^3} \b[\p_r(\eta^2 r^{a_2-1})- \p^2_r(\eta^2 r^{a_2})\b] |\p_{rz}^2 h_{\th}|^2 dx\\\no
&\q&+2\int_{\mbR^3} \eta^2 r^{a_2-2} |\p_{rz}^2 h_{\th}|^2 dx
-\int_{\mbR^3} \eta^2 r^{a_2} \p_{rzz}^3 h_{\theta} (u_r\p_r+ u_z\p_z)\p_r h_{\th}dx\\\no
& \q &- \int_{\mbR^3} \eta^2 r^{a_2} \p_{rzz}^3 h_{\theta}(\p_r u_r\p_r + \p_r u_z\p_z) h_{\theta}dx
- \int_{\mbR^3} \eta^2 r^{a_2} \p^3_{rzz} h_{\th} \b(\f {u_r}{r}- \f {\p_r u_r}{r}\b)\f{h_{\th}}{r} dx\\\no
&\q& +\int_{\mbR^3} \eta^2 r^{a_2} \p^3_{rzz} h_{\th} \f{u_r}{r} \p_r h_{\th} dx
+2 \int_{\mbR^3} \eta^2 r^{a_2-2} \p^3_{rzz} h_{\th} \f{h_{\th}}{r} dx\\\no
&\q&- \int_{\mbR^3} \eta^2 r^{a_2} \p_{rzz}^3 h_{\th} \p_r (\p_z g_r - \p_r g_z) dx
\\\no
& \eqqcolon & - \int_{\mbR^3} \eta^2 r^{a_2} |\na \p_{rz}^2 h_{\th}|^2 dx + \sum_{i=1}^2 D_{2i} + \sum_{j=1}^6 E_{2j}.
\ee
Using previous estimates, we can bound these terms as follows.
\be\no
|D_{21}|&\leq& C \int_{r\geq r_0} r^{a_2-2} |\p_{rz}^2 h_{\th}|^2 dx\leq C \int_{\mbR^3} r^{1+\ga} | \p_{rz}^2 h_{\th}|^2 dx,\q\text{if } a_2\leq 3+\ga,\\\no
|D_{22}|&\leq& C \int_{\mbR^3} \eta^2 r^{a_2-2} |\p^2_{rz} h_{\th}|^2 dx \leq C \int_{\mbR^3} r^{1+\ga} |\p_{rz}^2 h_{\th}|^2 dx,\q \text{if }a_2\leq 3+\ga,
\\\no
|E_{21}|
&\leq&\f1{16}\int_{\mbR^3} \eta^2 r^{a_2}|\na\p_{rz}^2 h_{\th}|^2 dx+  \int_{\mbR^3} \eta^2 r^{a_2}(|u_r|+|u_z|)^2  |\p_{rz}^2 h_{\th}|^2 dx\\\no
&\leq& \f{1}{16} \int_{\mbR^3} \eta^2 r^{a_2} |\na \p^2_{rz} h_{\th}|^2 dx + C \int_{\mbR^3} r^{1+\ga} |\p^2_{rz} h_{\th}|^2 dx,\q\text{if } a_2\leq 1+\ga+2\de,
\\\no
|E_{22}|
&\leq&\f1{16}\int_{\mbR^3} \eta^2 r^{a_2}|\na\p_{rz}^2 h_{\th}|^2 dx+  \int_{\mbR^3} \eta^2 r^{a_2}(|\na u_r|+|\na u_z|)^2 |\na h_{\th}|^2 dx\\\no
&\leq& \f{1}{16} \int_{\mbR^3} \eta^2 r^{a_2} |\na \p^2_{rz} h_{\th}|^2 dx + C \|\na h_{\th}\|_{L^2}^2,\q\text{if }a_2\leq 2(1+\ga),
\ee\be\no
|E_{23}|
&\leq&\f1{16}\int_{\mbR^3} \eta^2 r^{a_2}|\na\p_{rz}^2 h_{\th}|^2 dx+  \int_{\mbR^3} \eta^2 r^{a_2}\b(\f {|u_r|^2}{r^2} + |\na u_r|^2\b)|\f{h_{\theta}}{r}|^2 dx\\\no
&\leq& \f1{16}\int_{\mbR^3} \eta^2 r^{a_2}|\na\p_{rz}^2 h_{\th}|^2 dx+ C\int_{\mbR^3} \f {h_{\th}^2}{r^2},\q\text{if } a_2\leq 2\min\{1+\ga, 1+\de\},
\\\no
|E_{24}|
&\leq&\f1{16}\int_{\mbR^3} \eta^2 r^{a_2}|\na\p_{rz}^2 h_{\th}|^2 dx + C \int_{\mbR^3} \eta^2 r^{a_2} \f{u_r^2}{r^2} |\na h_{\th}|^2 dx,\q\text{if } a_2\le 2(1+\de),
\\\no
|E_{25}|
&\leq& \f1{16}\int_{\mbR^3} \eta^2 r^{a_2}|\na\p_{rz}^2 h_{\th}|^2 dx + C \int_{\mbR^3} \eta^2 r^{a_2-4} \f {h_{\th}^2}{r^2} dx\\\no
&\le& \f1{16}\int_{\mbR^3} \eta^2 r^{a_2}|\na\p_{rz}^2 h_{\th}|^2 dx + C \int_{\mbR^3} \f {h_{\th}^2}{r^2} dx,\q\text{if } a_2\le 4,\\\no
| E_{26}|
& \leq & \f1{16}\int_{\mbR^3} \eta^2 r^{a_2}|\na\p_{rz}^2 h_{\th}|^2 dx +  C \int_{\mbR^3} \eta^2 r^{a_2} |\p_r (\p_z g_r - \p_r g_z)|^2 dx.
\ee

We infer that
\be\lab{htheta33}
\int_{\mbR^3} r^{1+2\de +\ga}|\na \p_{rz}^2 h_{\th}|^2 dx<\oo.
\ee
Similarly, there also holds
\be\lab{htheta34}
\int_{\mbR^3} r^{1+2\de +\ga} |\na\p_{zz}^2 h_{\th}|^2 dx<\oo.
\ee
Therefore
\be\no
\int_{\mbR^3} r^{1+2\de +\ga}(|\na \p_z \om^h_r|^2+ |\na \p_z \om^h_z|^2) dx<\oo.
\ee
\end{proof}

By \eqref{uhdecay} and (\ref{na-u-rz}), choose $\delta=(\frac12)^-$ and $\gamma=(\frac14)^-$ in Lemma \ref{weighted estimate3}, and we are led to
\be\lab{htheta48}
\int_{\mbR^3} r^{(\f54)^-}(|\na \om^h_r|^2 +|\na \om^h_z|^2 ) dx<\oo, \\\lab{htheta49}
\int_{\mbR^3} r^{(\f 94)^-} (|\na \p_z \om^h_r|^2 + |\na \p_z \om^h_z|^2) dx<\oo,
\ee
and
\be\lab{htheta50}
|\om^h_r(r,z)| + |\om_z^h (r,z)|\leq C(M) r^{-(\f{35}{32})^-}.
\ee
Note that $\omega_r^h= -\p_z h_{\theta}$ and $\p_r h_{\theta}= \omega_z^h- \frac{h_{\theta}}{r}$, we also have (\ref{decay-na-h-th}). The proof of Theorem \ref{main12} is finished.

\section{Proof of Theorem \ref{main13}}\noindent

Note that the scaling technique by employing the \textit{Brezis-Gallouet inequality} can not be applied to get any decay along the axis, the only way we know is to do the weighted energy estimates with the weight $\rho = \sqrt{r^2 + z^2}$.

 The proof is organized as follows. In Step 1, some weighted energy estimates for $\Pi \= \f {h_{\th }}{r}$ were derived by using the weight $\rho=\sqrt{r^2+z^2}$ as in \cite{weng15}. In Step 2, we further infer that $|\Pi (r,z)|^2 \le \f{C}{r^2|z|^{(\f{31}{84}+\f{41}{84}\tau)^-}}$ and
\be\no
|h_{\theta}(r,z)| &\leq & C(M)(\rho+1)^{-(\f{31}{168}+\f{41}{168}\tau)^-},\q \forall (r,z)\in \mbR_+\times \mbR, \rho=\sqrt{r^2+z^2}.
\ee
In Step 3, where the special case $u_{\th} \equiv 0$ is considered, it will be shown that
\be\no
|\om_{\th}^u (r,z)| & \leq & C(M) (1+\rho)^{-(\f {13}{64} (1+\tau))^-}, \q \forall (r,z) \in \mbR_+\times \mbR, \rho = \sqrt{r^2 + z^2},
\ee
using similar arguments in the previous two steps on the equation of $\Om \= \f {\om_{\th}^u }r$. In Step 4, the decay of $(u_r, u_z)$ will be obtained from $\omega_{\theta}^u$ by the Biot-Savart law. Finally, we complete the proof by an iteration of $\tau$.

\begin{proof}[Proof of Theorem \ref{main13}]

{\bf Step 1.} We first notice that the quantity $\Pi \= \f{h_{\th}}{r}$ satisfies the following elliptic equation
\be\lab{Pi}
(u_r\p_r+ u_z\p_z)\Pi = \b(\p_r^2+\f{3}{r}\p_r+ \p_z^2\b)\Pi + \f 1r (\p_z g_r - \p_r g_z).
\ee
Following the argument developed in \cite{weng15}, we have the following weighted estimates for $\Pi$.

{\it Suppose that
\be\lab{urz100}
|u_r(r,z)| + |u_z(r,z)|\leq C (1+\rho)^{-\tau},\q \rho=\sqrt{r^2+z^2}
\ee
for some $\tau\in [0,1]$, and \emph{\bf{g}} satisfies (\ref{fg4}) and (\ref{fg5}), then we have
\be\lab{pi1}
&& \int_{\mbR^3} |\Pi(r,z)|^2 r dr dz < \oo,\\
&&\int_{\mbR^3} \rho^{1+\tau}|\nabla\Pi(r,z)|^2 rdr dz <\oo,\\\lab{pi2}
&&\int_{\mbR^3} \rho^{1+3\tau} |\nabla\p_z\Pi(r,z)|^2 rdr dz <\oo.
\ee
}

{\bf Step 2.} Derive the decay rate for $h_{\th}$.

By a priori estimates, we have
\be\lab{nahth}
&&\int_{\mbR^3} \f {|h_{\th}(r,z)|^2}{r^2}+ |\na h_{\th} (r,z)|^2 dx < \oo.
\ee

Combining the results in \eqref{naomh} and (\ref{pi1})-(\ref{nahth}), then
\be\lab{pi3}
&&\int_{\mbR^3} |\Pi(r,z)|^2 dx<\oo,\\\lab{pi4}
&&\int_{\mbR^3} (r^{2}+|z|^{1+\tau}) |\na\Pi(r,z)|^2 dx<\oo,\\\lab{pi5}
&&\int_{\mbR^3} (r^{3+\ga}+|z|^{1+3\tau}) |\na\p_z\Pi(r,z)|^2 dx<\oo,
\ee
where $\ga$ can be any constant less than $\f 14$. Fix $d>1$, then for each $n\in \mathbb{N}$,
\be\no
\int_{2^n}^{2^{n+1}} \int_d^{\oo} |\Pi(r,z)|^2 rdr dz<\oo.
\ee
By mean value theorem, there exists $z_n\in [2^n, 2^{n+1}]$ such that
\be\no
\int_d^{\oo} |\Pi(r,z_n)|^2 rdr \leq \f{C}{z_n}.
\ee
Then for any $z$, choose $z_n>z$ and
\be\no
\int_d^{\oo} |\Pi(r,z)|^2 r dr &=&\int_d^{\oo} |\Pi(r,z_n)|^2 r dr- 2\int_d^{\oo} \int_z^{z_n}\Pi(r,t)\p_t \Pi(r,t) r dr dt=: I_1 +I_2,\\\no
|I_2|&\leq&\b(\int_d^{\oo}\int_z^{z_n} |\Pi(r,t)|^2 r drdt\b)^{1/2}\b(\int_d^{\oo}\int_z^{z_n} |\p_t\Pi(r,t)|^2 r drdt\b)^{1/2}\\\no
&\leq& \f{C}{|z|^{\f12(1+\tau)}}.
\ee
Letting $z_n\to \oo$, then $I_1\to 0$ and
\be\lab{pi6}
\int_d^{\oo}|\Pi(r,z)|^2 rdr\leq \f{C}{|z|^{\f12(1+\tau)}}.
\ee

Similarly, one can find $z_n\in [2^n, 2^{n+1}]$ such that
\be\no
\int_d^{\oo} |\na\Pi(r,z_n)|^2 rdr &\leq& \f{C}{z_n^2},\\\no
\int_d^{\oo} |\na\Pi(r,z)|^2 r dr&=&\int_d^{\oo} |\na\Pi(r,z_n)|^2 r dr- 2\int_d^{\oo}\int_z^{z_n} \na\Pi(r,t)\cdot \p_t\na\Pi(r,t) r dr dt\\\no
& \eqqcolon & J_1+J_2,\\\no
|J_2|&\leq&\b(\int_d^{\oo}\int_z^{z_n} |\na\Pi(r,t)|^2 r drdt\b)^{1/2}\b(\int_d^{\oo}\int_z^{z_n} |\p_t\na\Pi(r,t)|^2 r drdt\b)^{1/2}\\\no
&\leq&\begin{cases}
\f{C}{d^{\f{5+\ga}{2}}},\\
\f{C}{|z|^{1+2\tau}}.
\end{cases}
\ee
Letting $n\to \oo$, $J_1\to 0$. Recall that $\ga=(\f 14)^-$, then $J_2\leq \min\b\{\f{C}{d^{\f{5+\ga}{2}}}, \f{C}{|z|^{1+2\tau}}\b\}$ and by interpolation
\be\no
\int_d^{\oo} |\na\Pi(r,z)|^2 r dr\leq \b(\f{C}{d^{(\f{21}{8})^-}}\b)^{\f{16}{21}} \b(\f{C}{|z|^{1+2\tau}}\b)^{\f{5}{21}}\leq \f{C}{d^2 |z|^{(\f5{21}(1+2\tau))^-}}.
\ee

Finally,
\be\no
|\Pi(d,z)|^2 &=&\f{1}{r_1-d}\int_d^{r_1} |\Pi(r,z)|^2 dr+ (|\Pi(r,z)|^2-\f{1}{r_1-d}\int_d^{r_1}|\Pi(r,z)|^2 dr)\\\no
& \eqqcolon & H_1 +H_2,
\ee
\be\no 
|H_2|&=&\b||\Pi(d,z)|^2-|\Pi(d_*,z)|^2\b| \leq 2\int_d^{r_1} |\Pi(r,z)\p_r\Pi(r,z)|dr\\\no
&\leq&\f{C}{d}\b(\int_d^{\oo}|\Pi(r,z)|^2 rdr\b)^{1/2}\b(\int_d^{\oo}|\na \Pi(r,z)|^2 r dr\b)^{1/2}\\\no
&\leq&\f{C}{d}\b(\f{C}{|z|^{\f12(1+\tau)}}\b)^{1/2}\b(\f{C}{d^2 |z|^{(\f 5{21} (1 + 2 \tau))^-}}\b)^{\f 12}\leq \f{C}{d^2|z|^{(\f{31}{84}+\f{41}{84}\tau)^-}},
\ee
It follows from \eqref{pi-z} that $\displaystyle\lim_{r_1\to \infty} H_1=0$ and
\be\lab{pi-z}
|h_{\theta}(d,z)|\leq C(M) |z|^{-(\f{31}{168}+\f{41}{168}\tau)^-}.
\ee

Together with previous decay results on $r$, we have
\be\lab{htheta-decay}
|h_{\theta}(r,z)| &\leq & C(M)(\rho+1)^{-(\f{31}{168}+\f{41}{168}\tau)^-},\q \forall (r,z)\in \mbR_+\times \mbR, \rho=\sqrt{r^2+z^2}.
\ee
In general case where $u_{\th} \neq 0$, we only have $\tau=0$, which implies that $|h_{\theta}| \le C(M) (\rho +1)^{- \f {31}{168}}$. We remark here that we do not obtain any decay rate for $\Pi$, since the weighted estimate for $\na \p_r \Pi$ is not available.

{\bf Step 3.} For the special case where the flow has zero swirl $u_{\th}\equiv 0$, we can further derive the decay rates of ${\bf u}$.

In this case, one has
\be\no
(u_r\p_r + u_z\p_z )\om_{\th}^u -\f{u_r}{r} \om^u_{\th} + \f1r\p_z (h_{\theta}^2) =\b(\p_r^2+\f1r\p_r+\p_z^2-\f1{r^2}\b) \om^u_{\th} - (\p_z f_r - \p_r f_z) ,
\ee
and
\be\no
(u_r\p_r+ u_z\p_z)\Om +\p_z \Pi^2= \b(\p_r^2+\f{3}{r}\p_r + \p_z^2\b)\Om + \f 1r (\p_z f_r -\p_r f_z) ,
\ee
for $\Om \= \f{\om^u_{\th}(r,z)}{r}$. We run the same argument as in \cite{weng15} for $\Om$ to derive the estimate $\int_{\mathbb{R}^3} \rho^{d_1}|\nabla \Omega|^2 dx$, there is an additional term
\be\no
\b|\int_{\mbR^3} \phi^2\rho^{d_1}\Om \p_z \Pi^2 dx\b| & \le & 2 \int_{\mbR^3} \phi^2 \rho^{d_1} |\Pi| |\Om| |\p_z \Pi| dx\\\no
 & \le & C \|\Pi \|_{L^{\oo}} \| \Om\|_{L^2} (\int_{\mbR^3} \rho^{2 d_1} |\p_z \Pi|^2 dx)^{1/2} \\\no
 & < & \oo,\q \text{if } d_1 \le \f {1+ \tau}{2},
\ee
where $\phi$ is a cut-off function satisfying $\phi \in C^{\oo}_0 (\mbR^3)$, $\phi = \phi (\rho)$, $0 \le \phi \le 1$, $\phi(\rho) =1$ on $2 \rho_0 \le \rho \le \rho_1$, $\phi (\rho) =0$ on $\rho \le \rho_0$ or $\rho \ge 2 \rho_1$, such that $|\na \phi| \le \f {C}{\rho_0}$ on $\rho_0 < \rho < 2 \rho_0$ and $|\na \phi| \le \f {C}{\rho_1}$ on $\rho_1 < \rho < 2 \rho_1$, and $\na \phi =0$ elsewhere.
Thus
\be\lab{Om1}
\int_{\mbR^3} \rho^{\f {1 + \tau}2} |\na \Om(r,z)|^2 dx<\oo.
\ee

To derive the estimate of $\na \p_z\Om$, we see that $\p_z\Om$ satisfies
\be\lab{pz-Om}
\p_z[(u_r\p_r + u_z\p_z) \Om + \p_z \Pi^2] = \b(\p_r^2 +\f{3}{r}\p_r+ \p_z^2\b)\p_z \Om + \p_z (\f 1r (\p_z f_r - \p_r f_z)).
\ee
Multiplying (\ref{pz-Om}) by $\phi^2 \rho^{d_2} \p_z \Om$ and integrating over $\mbR^3$, then we get
\be\no
0 & = & \int_{\mbR^3} \phi^2 \rho^{d_2} |\na \p_z \Om|^2 dx - \f 12 \int_{\mbR^3} [\p_r^2 (\phi^2 \rho^{d_2}) + \p_z^2 (\phi^2 \rho^{d_2})] |\p_z \Om|^2 dx  \\\no
&&+ \f 12 \int_{\mbR^3} \p_r (\phi^2 \rho^{d_2}) \f {|\p_z \Om|^2} {r} dx - \int_{\mbR^3} \p_z (\phi^2 \rho^{d_2}) \p_z \Om (u_r \p_r + u_z \p_z) \Om dx \\\no
&&- \int_{\mbR^3} \phi^2 \rho^{d_2} \p_z^2 \Om (u_r \p_r + u_z \p_z) \Om dx
 - \int_{\mbR^3} \p_z(\phi^2 \rho^{d_2}) \p_z \Om \cdot \p_z \Pi^2 dx\\\no
&& - \int_{\mbR^3} \phi^2 \rho^{d_2} \p_z^2 \Om \cdot \p_z \Pi^2 dx
 + \int_{\mbR^3} \p_z (\phi^2 \rho^{d_2}) \p_z \Om   (\f {(\p_z f_r - \p_r f_z)}{r}) dx \\\no
 && + \int_{\mbR^3} \phi^2 \rho^{d_2} \p_z^2 \Om \f  {(\p_z f_r - \p_r f_z)}{r}dx\\\no
& \eqqcolon & \int_{\mbR^3} \phi^2 \rho^{d_2} |\na \p_z \Om|^2 dx + \sum_{i=1}^6 K_i + \sum_{j=1}^2 Q_j.
\ee
These terms can be bounded as follows.
\be\no
|K_1| + |K_2| & \le & C \int_{\rho \ge \rho_0} \rho^{d_2-2} |\na \Om|^2 dx \le \oo, \q \text{if } d_2 \le \f {5 + \tau}{2},
\\\no
|K_3| & \le & C \int_{\rho \ge \rho_0} \rho^{d_2-1} |(u_r, u_z)| |\na \Om|^2 dx \le \oo, \q \text{if } d_2 \le \f 32(1+ \tau),\
\\\no
|K_4| & \le & \f 18 \int_{\mbR^3} \phi^2 \rho^{d_2} |\na \p_z \Om|^2 dx + C \int_{\mbR^3} \phi^2 \rho^{d_2} |(u_r, u_z)|^2 |\na \Om|^2 dx, \q \text{if } d_2 \le \f 12 + \f {5\tau} 2,
\\\no
|K_5|&\le& C \int_{\mbR^3} \phi \rho^{d_2-1} |\na\Om||\Pi||\p_z\Pi| dx\\\no
 &\leq& C \int_{\mbR^3} \rho^{\f {1 + \tau} 2} |\na\Om|^2 dx + C \int_{\mbR^3} \rho^{2(d_2-1)-\f {1+ \tau} 2} |\Pi|^2 |\p_z\Pi|^2 dx\\\no
&\leq& C \int_{\mbR^3} \rho^{1+\tau} |\na\Om|^2 dx + C \int_{\mbR^3} \rho^{1+\tau}|\p_z\Pi|^2 dx,\q\text{if } d_2\leq \f 74 + \f 34 \tau,\\\no
|K_6|&\leq& \f{1}{8} \int_{\mbR^3} \phi^2 \rho^{d_2}|\p_z^2\Om|^2 dx+ C \int_{\mbR^3}\phi^2 \rho^{d_2} \Pi^2 |\p_z \Pi|^2 dx\\\no
&\leq& \f{1}{8} \int_{\mbR^3} \phi^2 \rho^{d_2}|\p_z^2\Om|^2 dx+ C \int_{\mbR^3}\rho^{1+\tau}|\p_z \Pi|^2 dx, \q\text{if } d_2\leq 1 + \tau,
\ee
and
\be\no
| Q_1| & \le & C \int_{\rho \ge \rho_0} \rho^{2 (d_2 -1)- (1+\tau)} |\f {\p_z f_r - \p_r f_z}r|^2 dx + C \int_{\rho \ge \rho_0} \rho^{\f {1+ \tau} 2} |\na \Om|^2 dx, \q \text{if } d_2 \le \f {4+\tau}{2},\\\no
| Q_2 | & \le & \f 18 \int_{\mbR^3} \phi^2 \rho^{d_2} |\na \p_z \Om|^2 dx + C \int_{\mbR^3} \phi^2 \rho^{d_2}  |\f {\p_z f_r - \p_r f_z}r|^2 dx.
\ee
Letting $\rho_1\to \infty$, we obtain
\be\lab{Om2}
\int_{\mbR^3} \rho^{1+\tau} |\na \p_z \Om|^2 dx<\oo.
\ee

Combining the weighted energy estimate in Lemma \ref{weighted estimate1} and \eqref{Om1}-\eqref{Om2}, we find
\be\no
&&\int_{\mbR^3} r^2 |\Om|^2 dx<\oo,\\\no
&&\int_{\mbR^3} (r^{3+\de} +|z|^{ \f {1 + \tau}2 }) |\na \Om|^2 dx<\oo,\\\no
&&\int_{\mbR^3} (r^{3+3\de}+|z|^{1+\tau}) |\na \p_z \Om|^2 dx <\oo.
\ee

Same as {\bf Step 2}, fix any $d >1$, we have
\be\no
\int_{d}^{\oo} |\Om (r,z)|^2 r dr \le \f {C}{ d |z|^{\f {1+ \tau}{4}}}.
\ee
For each $n \in \mathbb{N}$,
\be\no
\int_{2^n}^{2^{n+1}} \int_d^{\oo} |z|^{\f {1+\tau}{2}} |\na \Om(r,z)|^2 r dr dz < \oo.
\ee
By mean value theorem, one can find $z_n \in [2^n, 2^{n+1}]$ such that
\be\no
\int_d^{\oo} |\na \Om (r, z_n)|^2 r dr & \le & \f {C}{|z_n|^{3/2}},\\\no
\int_d^{\oo} |\na \Om (r,z)|^2 r dr & =& \int_d^{\oo} |\na \Om (r,z_n)|^2 r dr - 2 \int_d^{\oo} \int_z^{z_n} \na \Om (r,t) \cdot \p_t \na \Om(r,t) r dr dt \\\no
& \=& L_1 + L_2,
\\\no
|L_2|&\leq&\b(\int_d^{\oo}\int_z^{z_n} |\na\Pi(r,t)|^2 r drdt\b)^{1/2}\b(\int_d^{\oo}\int_z^{z_n} |\p_t\na\Pi(r,t)|^2 r drdt\b)^{1/2}\\\no
&\leq&\begin{cases}
\f{C}{d^{3 + 2 \delta}},\\
\f{C}{|z|^{\f 34 (1 +\tau)}}.
\end{cases}
\ee
Letting $n \to \oo$, $L_1 \to 0$. Take $\delta = (\f 12)^-$ and $L_2 \le \textrm{min} \{\f{C}{d^{3 + 2 \delta}} ,\f{C}{|z|^{\f 34 (1 +\tau)}} \}$,
\be\no
\int_d^{\oo} |\na \Om (r,z)|^2 r dr \le \b(\f C{d^{4-}}\b)^{\f 14} \b( \f C {|z|^{\f 34 (1+\tau)}}\b)^{\f 34} \le \f C{d |z|^{\f 9{16} (1+\tau)^-}}.
\ee
Hence,
\be\no
|\Om(d,z)|^2 &=& \f 1{r_1 - d} \int_d^{r_1} |\Om(r,z)|^2 dr + \b(|\Om(d,z)|^2 - \f 1{r_1 - d} \int_d^{r_1} |\Om(r,z)|^2 dr \b) \\\no
& \eqqcolon & G_1 + G_2.
\ee
\be\no 
|G_2| & =& \b||\Om(d,z)|^2 - |\Om(d_*,z)|^2 \b| \\\no
& \le & 2 \int_d^{r_1} |\Om(r,z) \cdot \p_r \Om(r,z)| dr\\\no
&\le & \f Cd \b(\int_d^{\oo} |\Om(r,z)|^2 r dr \b)^{1/2} \b(\int_d^{\oo} |\na \Om(r,z)|^2 r dr \b)^{1/2}\\\no
& \le & \f Cd \b(\f C{d |z|^{\f 14 (1+\tau)^-}} \b)^{1/2} \b(\f C{d |z|^{\f 9 {16} (1+ \tau)^-}} \b)^{1/2} \le \f {C}{d^2 |z|^{\f {13}{32} (1+\tau)^-}}.
\ee
Therefore we have
\be\no
|\om_{\th}^u (r,z)| & \le & C(M) (1+\rho)^{-(\f {13}{64} (1+\tau))^-}, \q \forall (r,z) \in \mbR_+\times \mbR, \rho = \sqrt{r^2 + z^2}.
\ee

{\bf Step 4.} Decay rate of ${\bf u}$.

Following the argument developed in \cite{weng15}, fix any $\textbf{x} \in \mbR^3 \setminus\{0\}$, define a cut-off function $\varphi \in C_0^{\oo} (\mbR^3)$ satisfying $\varphi({\bf y}) \equiv 1$ for $\forall \, {\bf y} \in B_{\rho/4}({\bf x})$ and $\varphi ({\bf y}) \equiv 0$ for $\forall \,{\bf y} \notin B_{\rho/2}({\bf x})$, where $\rho = |{\bf x}|$ and $|\na \varphi ({\bf y})| \le \f C{|y|}$, $|\na^2 \varphi ({\bf y})| \le \f C{|y|^2}$ for ${\bf y} \in D:= B_{\rho/2}({\bf x}) \setminus B_{\rho/4} ({\bf x})$. For ${\bf u} ({\bf x}) = u_r {\bf e}_r + u_z {\bf e}_z $, since $\textrm{curl } {\bf u} = \om_{\th}^u {\bf e}_{\th}$, then
\be\no
{\bf u}({\bf x}) & = & - \int_{\mbR^3} \na_{\bf y} \Ga ({\bf x},{\bf y})\times (\om_{\phi}^u({\bf y}) \varphi({\bf y}) {\bf e}_{\phi}) d{\bf y}
- \int_{\mbR^3} \Ga ({\bf x},{\bf y}) ( \na_{\bf y}\varphi({\bf y}) \times {\bf e}_{\phi} ) \om_{\phi}^u({\bf y}) d{\bf y} \\\no
&& \int_{\mbR^3} \Ga ({\bf x},{\bf y}) (\De_{\bf y} \varphi) ({\bf y}) {\bf u} ({\bf y}) d{\bf y}
+ 2 \int_{\mbR^3} (\na_{\bf y} \Ga )({\bf x},{\bf y}) \cdot (\na_{\bf y} \varphi ) ({\bf y}) {\bf u} ({\bf y}) d{\bf y}
\\\no
& \= & \sum_{i=1}^4 R_i
\ee
According to \cite{weng15}, we have $|R_i| \le \f C{\rho^{1/2}}, \q i=2,3,4$. For the estimate of $R_1$, fix a $d \in (0, \f {\rho}{2})$, which will be determined later, then
\be\no
|R_1| & \le &  \sup \limits_{{\bf y} \in B_d ({\bf x})} |\om_{\phi}^u({\bf y})| \int_{B_d ({\bf x})} |\na \Ga ({\bf x}-{\bf y})| d{\bf y}\\\no
&& + \b(\int_{B_{\rho/2} ({\bf x}) \setminus B_d ({\bf x})}|\na \Ga ({\bf x}-{\bf y})|^2 d{\bf y}  \b)^{\f 12} \b(\int_{B_{\rho/2} ({\bf x}) \setminus B_d ({\bf x})}|\om_{\phi}^u({\bf y})|^2 d{\bf y} \b)^{\f 12}\\\no
& \le & C\rho^{-(\f {13}{64} (1+\tau))^-}d + C d^{-\f 12}.
\ee
By choosing $d= \rho^{\f {13}{96} (1+ \tau)^-}$, we obtain the bound for $|R_1| \le \f C{\rho^{\f {13}{192} (1+\tau)^-}}$. Hence we have
\be
|(u_r, u_z) (r,z)| \le C(M) (1+ \rho)^{- \f{13}{192}(1+\tau)^-}.
\ee

{\bf Step 5.} Iteration.

At the beginning, we have $\tau=0$ in \eqref{urz100}, then by using the arguments developed in {\bf Step 1} to {\bf Step 4}, we have a new $\tau$ in \eqref{urz100}, which will be denoted by $\tau_1=(\f{13}{192})^-$. Running a second iteration of these three steps, we get a new $\tau_2= \tau_1+ \f {13}{192}\tau_1$, and after $n$ iteration, we get
\be\no
\tau_n= \tau_{1}+ \f {13}{192}\tau_{n-1}= \tau_1 \sum_{i=0}^{n-1}(\f {13}{192})^i.
\ee
Let $n\to \oo$, $\tau_n\to (\f {13}{179})^-$ as $n\to \oo$. In a word, we infer the following decay rates
\be\lab{urz300}
|u_r(r,z)| +|u_z(r,z)|&\leq& C(M) (\rho+1)^{-(\f {13}{179})^-},\\\lab{omgea300}
|\om_{\th}(r,z)| &\leq& C(M) (\rho+1)^{-(\f {39}{179})^-},\\\lab{htheta30}
|h_{\th}(r,z)|&\leq & C(M) (\rho + 1)^{-(\f {3041}{15036})^-}.
\ee
\end{proof}

{\bf Acknowledgement.} Weng is partially supported by National Natural Science Foundation of China 11701431, 11971307, 12071359.  We sincerely thank the anonymous reviewers for their careful reading and thoughtful suggestions that help to improve this paper.

\end{document}